\documentclass[10.5pt]{article}
\usepackage{amsmath}
\usepackage{amssymb}
\usepackage{amscd}
\usepackage{xspace}
\usepackage{verbatim}
\renewcommand{\theequation}{\thesection.\arabic{equation}
}
\newcommand{\mysection}[1]{\section{#1}\setcounter{equation}{0}}
\title{\bf Initial trace of solutions of semilinear heat equation with absorption}
\author{{\bf Waad Al Sayed\footnote{Email: waadalsayed@hotmail.com}}\\
{\small College of Sciences and Humanities, Fahed Bin Sultan University,  Tabuk,  Saudi arabia}\\[2mm]
 {\bf Laurent V\'eron\footnote{Email: veronl@univ-tours.fr}}\\
{\small Department of Mathematics, Universit\'e Fran\c{c}ois Rabelais,  Tours,  France}
}
\date{}
\begin{document}
\maketitle
{\small{\bf Abstract} We study the initial trace problem for positive solutions of semilinear heat equations with strong absorption. We show that in general this initial trace is an outer regular Borel measure. We emphasize in  particular the case where $u$ satisfies (E) $\partial_t u-\Delta u+t^\alpha |u|^{q-1}u=0$, with $q>1$ and $\alpha>-1$ and prove that in the subcritical case $1<q<q_{\alpha,N}:=1+2(1+\alpha)/N$ the initial trace establishes a one to one correspondence between the set of outer regular Borel in measures $\mathbb R^N$ and the set of positive solutions of (E) in $\mathbb R^N\times\mathbb R_+$.
}\smallskip

\noindent Keywords:  Heat equation;  Singularities; Harnack inequalities; Radon measures; Borel measures; barrier; initial trace.\smallskip

\noindent MSC2010: 35K58; 35K15
\vspace{1mm}
\hspace{.05in}
\newcommand{\txt}[1]{\;\text{ #1 }\;}
\newcommand{\tbf}{\textbf}
\newcommand{\tit}{\textit}
\newcommand{\tsc}{\textsc}
\newcommand{\trm}{\textrm}
\newcommand{\mbf}{\mathbf}
\newcommand{\mrm}{\mathrm}
\newcommand{\bsym}{\boldsymbol}
\newcommand{\scs}{\scriptstyle}
\newcommand{\sss}{\scriptscriptstyle}
\newcommand{\txts}{\textstyle}
\newcommand{\dsps}{\displaystyle}
\newcommand{\fnz}{\footnotesize}
\newcommand{\scz}{\scriptsize}
\newcommand{\be}{
\begin{equation}
}
\newcommand{\bel}[1]{\begin{equation}\label{#1}}
\newcommand{\ee}{\end{equation}}
\newcommand{\eqnl}[2]{\begin{equation}\label{#1}{#2}\end{equation}}
\newtheorem{subn}{\name}
\renewcommand{\thesubn}{}
\newcommand{\bsn}[1]{\def\name{#1}
\begin{subn}}
\newcommand{\esn}{
\end{subn}}
\newtheorem{sub}{\name}[section]
\newcommand{\dn}[1]{\def\name{#1}}   
\newcommand{\bs}{
\begin{sub}}
\newcommand{\es}{
\end{sub}}
\newcommand{\bsl}[1]{
\begin{sub}\label{#1}}
\newcommand{\bth}[1]{\def\name{Theorem}
\begin{sub}\label{t:#1}}
\newcommand{\blemma}[1]{\def\name{Lemma}
\begin{sub}\label{l:#1}}
\newcommand{\bcor}[1]{\def\name{Corollary}
\begin{sub}\label{c:#1}}
\newcommand{\bdef}[1]{\def\name{Definition}
\begin{sub}\label{d:#1}}
\newcommand{\bprop}[1]{\def\name{Proposition}
\begin{sub}\label{p:#1}}
\newcommand{\R}{\eqref}
\newcommand{\rth}[1]{Theorem~\ref{t:#1}}
\newcommand{\rlemma}[1]{Lemma~\ref{l:#1}}
\newcommand{\rcor}[1]{Corollary~\ref{c:#1}}
\newcommand{\rdef}[1]{Definition~\ref{d:#1}}
\newcommand{\rprop}[1]{Proposition~\ref{p:#1}}
\newcommand{\BA}{
\begin{array}}
\newcommand{\EA}{
\end{array}}
\newcommand{\BAN}{\renewcommand{\arraystretch}{1.2}
\setlength{\arraycolsep}{2pt}
\begin{array}}
\newcommand{\BAV}[2]{\renewcommand{\arraystretch}{#1}
\setlength{\arraycolsep}{#2}
\begin{array}}
\newcommand{\BSA}{
\begin{subarray}}
\newcommand{\ESA}{
\end{subarray}}
\newcommand{\BAL}{
\begin{aligned}}
\newcommand{\EAL}{
\end{aligned}}
\newcommand{\BALG}{
\begin{alignat}}
\newcommand{\EALG}{
\end{alignat}}
\newcommand{\BALGN}{
\begin{alignat*}}
\newcommand{\EALGN}{
\end{alignat*}}
\newcommand{\note}[1]{\textit{#1.}\hspace{2mm}}
\newcommand{\Proof}{\note{Proof}}
\newcommand{\qeda}{\hspace{10mm}\hfill $\square$}
\newcommand{\qed}{\\
${}$ \hfill $\square$}
\newcommand{\modin}{$\,$\\
[-4mm] \indent}
\newcommand{\forevery}{\quad \forall}
\newcommand{\set}[1]{\{#1\}}
\newcommand{\setdef}[2]{\{\,#1:\,#2\,\}}
\newcommand{\setm}[2]{\{\,#1\mid #2\,\}}
\newcommand{\lra}{\longrightarrow}
\newcommand{\lla}{\longleftarrow}
\newcommand{\llra}{\longleftrightarrow}
\newcommand{\Lra}{\Longrightarrow}
\newcommand{\Lla}{\Longleftarrow}
\newcommand{\Llra}{\Longleftrightarrow}
\newcommand{\warrow}{\rightharpoonup}
\newcommand{
\paran}[1]{\left (#1 \right )}
\newcommand{\sqbr}[1]{\left [#1 \right ]}
\newcommand{\curlybr}[1]{\left \{#1 \right \}}
\newcommand{\abs}[1]{\left |#1\right |}
\newcommand{\norm}[1]{\left \|#1\right \|}
\newcommand{\paranb}[1]{\big (#1 \big )}
\newcommand{\lsqbrb}[1]{\big [#1 \big ]}
\newcommand{\lcurlybrb}[1]{\big \{#1 \big \}}
\newcommand{\absb}[1]{\big |#1\big |}
\newcommand{\normb}[1]{\big \|#1\big \|}
\newcommand{
\paranB}[1]{\Big (#1 \Big )}
\newcommand{\absB}[1]{\Big |#1\Big |}
\newcommand{\Remark}{\note{Remark}}
\newcommand{\normB}[1]{\Big \|#1\Big \|}
\newcommand{\thkl}{\rule[-.5mm]{.3mm}{3mm}}
\newcommand{\thknorm}[1]{\thkl #1 \thkl\,}
\newcommand{\trinorm}[1]{|\!|\!| #1 |\!|\!|\,}
\newcommand{\bang}[1]{\langle #1 \rangle}
\def\angb<#1>{\langle #1 \rangle}
\newcommand{\vstrut}[1]{\rule{0mm}{#1}}
\newcommand{\rec}[1]{\frac{1}{#1}}
\newcommand{\opname}[1]{\mbox{\rm #1}\,}
\newcommand{\supp}{\opname{supp}}
\newcommand{\dist}{\opname{dist}}
\newcommand{\myfrac}[2]{{\displaystyle \frac{#1}{#2} }}
\newcommand{\myint}[2]{{\displaystyle \int_{#1}^{#2}}}
\newcommand{\mysum}[2]{{\displaystyle \sum_{#1}^{#2}}}
\newcommand {\dint}{{\displaystyle \int\!\!\int}}
\newcommand{\q}{\quad}
\newcommand{\qq}{\qquad}
\newcommand{\hsp}[1]{\hspace{#1mm}}
\newcommand{\vsp}[1]{\vspace{#1mm}}
\newcommand{\ity}{\infty}
\newcommand{\prt}{\partial}
\newcommand{\sms}{\setminus}
\newcommand{\ems}{\emptyset}
\newcommand{\ti}{\times}
\newcommand{\pr}{^\prime}
\newcommand{\ppr}{^{\prime\prime}}
\newcommand{\tl}{\tilde}
\newcommand{\sbs}{\subset}
\newcommand{\sbeq}{\subseteq}
\newcommand{\nind}{\noindent}
\newcommand{\ind}{\indent}
\newcommand{\ovl}{\overline}
\newcommand{\unl}{\underline}
\newcommand{\nin}{\not\in}
\newcommand{\pfrac}[2]{\genfrac{(}{)}{}{}{#1}{#2}}
\newcommand{\dsp}{\displaystyle}
\newtheorem{Theo}{Theorem}[section]
\newtheorem{lema}{Lemma}[section]
\newtheorem{De}{Definition}[section]
\newtheorem{prop}{Proposition}[section]
\newtheorem{propriete}{Propri\'et\'e}[section]
\newtheorem{defi}{D\'efinition}[section]
\newtheorem{Proposition}{Proposition}[section]
\newtheorem{Lemme}{Lemme}[section]
\newtheorem{cor}{Corollary}[section]
\newtheorem{rem}{Remark}[section]
\def\ga{\alpha}     \def\gb{\beta}       \def\gg{\gamma}
\def\gc{\chi}       \def\gd{\delta}      \def\ge{\epsilon}
\def\gth{\theta}                         \def\vge{\varepsilon}
\def\gf{\phi}       \def\vgf{\varphi}    \def\gh{\eta}
\def\gi{\iota}      \def\gk{\kappa}      \def\gl{\lambda}
\def\gm{\mu}        \def\gn{\nu}         \def\gp{\pi}
\def\vgp{\varpi}    \def\gr{\rho}        \def\vgr{\varrho}
\def\gs{\sigma}     \def\vgs{\varsigma}  \def\gt{\tau}
\def\gu{\upsilon}   \def\gv{\vartheta}   \def\gw{\omega}
\def\gx{\xi}        \def\gy{\psi}        \def\gz{\zeta}
\def\Gg{\Gamma}     \def\Gd{\Delta}      \def\Gf{\Phi}
\def\Gth{\Theta}
\def\Gl{\Lambda}    \def\Gs{\Sigma}      \def\Gp{\Pi}
\def\Gw{\Omega}     \def\Gx{\Xi}         \def\Gy{\Psi}
\def\CS{{\mathcal S}}   \def\CM{{\mathcal M}}   \def\CN{{\mathcal N}}
\def\CR{{\mathcal R}}   \def\CO{{\mathcal O}}   \def\CP{{\mathcal P}}
\def\CA{{\mathcal A}}   \def\CB{{\mathcal B}}   \def\CC{{\mathcal C}}
\def\CD{{\mathcal D}}   \def\CE{{\mathcal E}}   \def\CF{{\mathcal F}}
\def\CG{{\mathcal G}}   \def\CH{{\mathcal H}}   \def\CI{{\mathcal I}}
\def\CJ{{\mathcal J}}   \def\CK{{\mathcal K}}   \def\CL{{\mathcal L}}
\def\CT{{\mathcal T}}   \def\CU{{\mathcal U}}   \def\CV{{\mathcal V}}
\def\CZ{{\mathcal Z}}   \def\CX{{\mathcal X}}   \def\CY{{\mathcal Y}}
\def\CW{{\mathcal W}} \def\CQ{{\mathcal Q}}
\def\BBA {\mathbb A}   \def\BBb {\mathbb B}    \def\BBC {\mathbb C}
\def\BBD {\mathbb D}   \def\BBE {\mathbb E}    \def\BBF {\mathbb F}
\def\BBG {\mathbb G}   \def\BBH {\mathbb H}    \def\BBI {\mathbb I}
\def\BBJ {\mathbb J}   \def\BBK {\mathbb K}    \def\BBL {\mathbb L}
\def\BBM {\mathbb M}   \def\BBN {\mathbb N}    \def\BBO {\mathbb O}
\def\BBP {\mathbb P}   \def\BBR {\mathbb R}    \def\BBS {\mathbb S}
\def\BBT {\mathbb T}   \def\BBU {\mathbb U}    \def\BBV {\mathbb V}
\def\BBW {\mathbb W}   \def\BBX {\mathbb X}    \def\BBY {\mathbb Y}
\def\BBZ {\mathbb Z}
\def\GTA {\mathfrak A}   \def\GTB {\mathfrak B}    \def\GTC {\mathfrak C}
\def\GTD {\mathfrak D}   \def\GTE {\mathfrak E}    \def\GTF {\mathfrak F}
\def\GTG {\mathfrak G}   \def\GTH {\mathfrak H}    \def\GTI {\mathfrak I}
\def\GTJ {\mathfrak J}   \def\GTK {\mathfrak K}    \def\GTL {\mathfrak L}
\def\GTM {\mathfrak M}   \def\GTN {\mathfrak N}    \def\GTO {\mathfrak O}
\def\GTP {\mathfrak P}   \def\GTR {\mathfrak R}    \def\GTS {\mathfrak S}
\def\GTT {\mathfrak T}   \def\GTU {\mathfrak U}    \def\GTV {\mathfrak V}
\def\GTW {\mathfrak W}   \def\GTX {\mathfrak X}    \def\GTY {\mathfrak Y}
\def\GTZ {\mathfrak Z}   \def\GTQ {\mathfrak Q}
\font\Sym= msam10
\def\SYM#1{\hbox{\Sym #1}}
\newcommand{\bdw}{\prt\Gw\xspace}

\mysection{Introduction}
\def\theequation{1.\arabic{equation}}\makeatother
In this paper we study the initial trace problem for positive solutions of 
\begin{equation}\label{0.0}
\partial _tu-\Delta u+ g(x,t,u)=0\mbox{ dans }\; Q^\Gw_T:=\Omega \times (0,T)
\end{equation}
where $\Omega$ is an open domain in $\mathbb R^N$, $g\in C(\Gw\ti\BBR_+\times \mathbb R)$ such that $g(x,t,.)$ is nondecreasing $\forall (x,t) \in \Gw\times \mathbb R$ and  $rg(x,t,r)\geq 0$ for all $(x,t,r)\in\Gw\ti\BBR_+\times \mathbb R$. Our first result establishes the existence of an initial trace. \smallskip


\noindent{ \bf Theorem A} {\it Assume $g$ satisfies the above conditions and that equation $(\ref{0.0})$ possesses a barrier at any $z\in\Gw$. If $u\in C^1(Q_T^\Gw)$ is a positive solution of (\ref{0.0}), it admits an initial trace which belongs to the class of outer regular positive Borel measure in $\Gw$.}\smallskip

The barrier assumption will be made precise later on in full generality. It is fulfilled if 
$g(x,t,r)=h(x)t^\ga\abs r^{q-1}r$ with $\ga>-1$, $q>1$ and $h\in L^\infty_{loc}(\Gw)$ satisfies 
$\inf\,{\rm \!ess}\,h>0$ for any compact subset $K\subset\Gw$, or if $g$ satisfies the Keller-Osserman condition, that is $g(x,t,r)\geq h(r)\geq 0$ where $h$ is nondecreasing and there exists $a$ such that
\begin{equation}\label{EE1}
\myint{a}{\infty}\frac{ds}{\sqrt{H(s)}}\qquad\text{where }H(s)=\myint{0}{s}h(t) dt.
\end{equation}

The initial trace of positive solutions of (\ref{0.0}) exists in the following sense: {\it there exists a relatively closed set $\CS\subset\Gw$ and a positive Radon measure $\gm$ on $\CR:=\Gw\setminus\CS$ with the following properties:\smallskip

\noindent (i) for any $x_0\in\CS$ and any $\ge>0$
 \begin{equation}\label{tr1}
\lim_{t\to 0}\myint{B_\ge(x_0)\cap\Gw}{}u(x,t)dx=\infty,
 \end{equation}
 \smallskip

\noindent (ii) for any $\gz\in C_c(\CR)$}
\begin{equation}\label{tr2}
\lim_{t\to 0}\myint{\Gw}{}u(x,t)\gz (x)dx=\myint{\Gw}{}\gz d\gm.
 \end{equation}

The couple $(\CS,\gm)$ is unique and characterizes a unique positive outer regular Borel measure $\gn$ on $\Gw$.  \smallskip

A similar notion of boundary trace has been introduced by Marcus and V\'eron \cite{MV3} in the study of  positive solutions of
\begin{equation}\label{EE}
-\Delta u+ g(x,u)=0\mbox{ in }\; \Gw.
\end{equation}
This notion in itself has turned out to be a very usefull tool for classifying the positive solutions of (\ref{EE}).\smallskip

In the second part we concentrate on the particular case of equation  
\begin{equation}\label{1}
\partial _tu-\Delta u+ t^{\alpha}\abs{u}^{q-1}u=0\mbox{ dans }\; Q^\Gw_T
\end{equation}
where $T>0$, $\alpha >-1 $ and $q>1$. Among the most useful tools, we point out the description of positive solutions with an isolated singularity 
at $(a,0)$ for $a\in\Gw$, whenever they exist: they are solutions $u$ of (\ref{1})   in $Q^\Gw_T$,  which belong to $C^{2,1}\left(Q^\Gw_T\right)\cap 
C\left(\Gw\times [0,T)\backslash \{(a,0)\}\right)$ and satisfy 
\begin{equation}\label{en0}
u(x,0)=0 \mbox{ in } \Gw \backslash \{a\}.
\end{equation}
When $\ga=0$, Brezis and Friedman prove in \cite{BF} that if $B_{2R}(a)\subset\Gw$,  then any such solution satisfies 
\begin{equation}\label{BFest}
u(x,t)\leq \frac{C\left(N,q,R\right)}{\left(\left|x-a\right|^2+t\right)^{\frac{1}{q-1}}}\quad\forall (x,t)\in B_R(a)\backslash\{0\}\ti [0,T].
\end{equation}
They also prove that if $1<q<q_N:=1+\frac{2}{N}$ and $k>0$ there exist singular solutions with initial data $u(.,0)=k\gd_a$, unique if $u$ vanishes on $\prt\Gw\times[0,T]$. In this range of exponents, Brezis, Peletier, Terman obtain in  \cite{BPT} the existence and uniqueness of a {\it very singular solution} of  (\ref{1}), always with $\ga=0$: it is a positive solution in $Q_\infty:=Q_\infty^{\BBR^N}$ under the form 
$$v_0(x,t)=t^{-1/(q-1)}V_0\left(\frac{x}{\sqrt t}\right),
$$
where $V_0>0$ is $C^2$ and satisfies 
\begin{equation}\label{sing}\BA {l}
-\Delta V_0 -\frac{1}{2}\eta.\nabla V_0 -\frac{1}{q-1} V_0+V_0^q=0 \quad\mbox{in } \mathbb R^N\\
\phantom{-,,;,--}\lim_{|\eta|\to\infty}\abs\eta^{\frac{2}{q-1}}V_0(\eta)=0.
\EA\end{equation}
Actually, Kamin and Peletier show that $v_0$ is the limit of the solutions $u_{k}$ of (\ref{1}) in $Q_\infty$ which satisfy 
$u(.,0)=k\gd_0$. The very singular singular solution plays a fundamental role in Marcus and V\'eron's description  \cite{MV1} of the initial trace of positive solutions of (\ref{1}) with $\ga=0$. In \cite{MV2}, Marcus et V\'eron study this equation when $\alpha\geq 0$ and $1<q<q_{\alpha,N}= 1+\frac{2(1+\alpha)}{N}$. They obtain the existence of a self-similar solution of (\ref{1}) in $Q_\infty$ under the form
$$v_\ga(x,t)=t^{-\frac{1+\alpha}{q-1}} V_\ga\left(\frac{x}{\sqrt t}\right),$$ 
which satisfies
$$\lim _{t\rightarrow 0}v_\ga(x,t)=0\;\;\;\forall x \not= 0$$
and
$$
\lim _{t\rightarrow 0}\myint{B_\ge}{}v_\ga(x,t)dx=\infty\;\;\;\forall \ge > 0.$$
The function $V_\ga$ is nonnegative and  verifies 
\begin{equation}\label{sing}
-\Delta V_\ga -\frac{1}{2}\eta.\nabla V_\ga -\frac{1+\ga}{q-1}  V_\ga+V_\ga^q=0 \mbox{ in } \mathbb R^N.
\end{equation}
Furthermore
\begin{equation}\label{sing+}V_\ga(\eta)=C\dsp|\eta|^{\frac{2(1+\alpha)}{q-1}-N}e^{\frac{-|\eta|^2}{4}}\left(1+o(1)\right) \mbox{ as }|\eta|\rightarrow \infty.
\end{equation}
If $1<q<q_{\alpha,N}$, they show that for every $k>0$ there exists a unique solution $u_{k\delta _a}$ of (\ref{1}) in $Q_\infty$ with initial data $k\delta _a$. Furthermore $\lim _{k\rightarrow \infty}u_{k\delta _a}=v_\ga$. Actually the limitation $\ga\geq 0$ can be relaxed to $\ga>-1$ has we will see it later on. Furthermore $Q_\infty$ can be replaced by $Q_\infty^\Gw$ provided $\prt\Gw$ is compact and smooth enough and $u_{k\delta _a}$ vanishes on $\prt\Gw\ti[0,\infty)$. \smallskip
 
 In this article we extend Brezis-Friedman removability result to equation (\ref{1}). We also extend Oswald's  classification of positive isolated singularities \cite{Osw}.  The starting point of our study is the following extension of estimate (\ref{BFest}) valid for any $\ga>-1$ and $q>1$.
 \begin{equation}\label{alpha-est}
u(x,t)\leq \frac{C\left(N,q,\ga,R\right)}{\left(\left|x-a\right|^2+t\right)^{\frac{1+\ga}{q-1}}}
\quad\forall (x,t)\in B_R\backslash\{a\}\ti [0,T].
 \end{equation}
 
The obstacle for obtaining such an estimate arises when $\ga>0$ and the absorption term $t^\ga u^q$ is degenerate near $t=0$. We overcome this difficulty by a delicate construction of 1-dim self-similar supersolutions. Thanks to this estimate, we obtain that the following classification result holds.\smallskip

\noindent{ \bf Theorem B} {\it Assume $1<q<q_{\ga,N}$ and $u\in C^1(Q_T^\Gw)\cap C(\Gw\ti [0,T]\backslash\{(a,0))\})$ is a solution of (\ref{1}) which vanishes on $\Gw\setminus\{a\}$ at $t=0$. Then\smallskip

\noindent (i) either there exists $k\geq 0$ such that $u(.,0)=k\gd_a$ and 
 \begin{equation}\label{weak}
u(x,t)\sim kE(x-a,t)\qquad\text{as }(x,t)\to (a,0),
 \end{equation}
where $E(x,t)=(4\pi t)^{-\frac{N}{2}}e^{-\frac{\abs x^2}{4t}}$,\smallskip

\noindent (ii) or  
 \begin{equation}\label{strong}
u(x,t)\sim v_\ga(x-a,t)\qquad\text{as }(x,t)\to (a,0).
 \end{equation}
}\\

In the supercritical case the following  removability statement holds.\smallskip

\noindent{ \bf Theorem C} {\it Assume $q\geq q_{\ga,N}$ and $u\in C^1(Q_T^\Gw)\cap C(\Gw\ti [0,T]\backslash\{(a,0))\})$ is a solution of (\ref{1}) which vanishes on $\Gw\setminus\{a\}$ at $t=0$. Then $u$ can be extended by continuity as a function in  $C(\Gw\ti [0,T])$}.\medskip

We prove that equation (\ref{1}) admits a barrier at any $z\in\Gw$. More precisely we construct a positive solution $w_{B_R}$  of (\ref{1}) in $Q^{B_R}_\infty$ which tends to $0$ locally uniformly in $B_R$ when $t\to 0$ and which blows-up uniformly on $\prt B_R\ti [\gt,\infty)$, for any $\gt>0$. Applying Theorem A, we infer that any positive solution admits an initial trace which is an outer regular Borel measure $\gn\approx(\CS,\gm)$. Using sharp parabolic Harnack inequality and a concentration principle, we prove the following result which is the key-stone  for analyzing the behaviour of $u$ on the set $\CS$.\smallskip   

 \smallskip

\noindent{ \bf Theorem D} {\it Assume $1<q< q_{\ga,N}$ and $u\in C^{2,1}(Q_T^\Gw)$ is a positive solution of (\ref{1}) with initial trace $(\CS,\gm)$. Then for any $a\in\CS$ there holds
\bel{sous}
u(x,t)\geq u_{\infty,a}(x,t)\qquad\forall (x,t)\in Q_T^\Gw
\ee
where $u_{\infty,a}=\lim_{k\to\infty}u_{k\gd_a}$,  $u_{k\gd_a}$being the solution of (\ref{1}) in $Q_T^\Gw$ with initial trace $k\gd_a$ and which vanishes on $\prt\Gw\ti [0,T]$. 
} \smallskip

It is important to notice that the behaviour of $u_{\infty,a}$ near $(a,0)$ is given by (\ref{strong}) and (\ref{sing+}).
Using (\ref{sous}), (\ref{alpha-est}) and sharp asymptotics of the function $V_\ga$, we are able to prove the following result which extends Theorem A.\smallskip 

\noindent{ \bf Theorem E} {\it Assume $1<q< q_{\ga,N}$ and $\Gw\subset\BBR^N$ is open with a $C^2$ compact boundary, eventually empty. Then for any couple $(\CS,\gm)$
where $\CS\subset\Gw$ is relatively closed and $\gm\in \mathfrak M_+(\Gw\setminus\CS)$, there exists  a maximal positive solution $\overline u$ and a minimal positive solution $\underline u$ of (\ref{1}), which belong to $C^{2,1}(Q^\Gw_T)\cap C(\overline\Gw\ti (0,T])$, satisfy (\ref{tr1}) and (\ref{tr2}) and vanish on $\prt\Gw\ti (0,T]$. If $\inf\{\abs{z-z'}:z\in \CS,z'\in\Gw^c\}>0$ and $\gm$ is bounded in a neighborhood of $\prt\Gw$, then $\overline u=\underline u$.}\medskip


\mysection{Initial trace}
\def\theequation{2.\arabic{equation}}\makeatother
In this section $\Gw\subset\BBR^N$ is an open domain, $ Q^{\Gw}_T=\Gw\ti (0,T)$, $\prt_\ell Q^{\Gw}_T=\Gw\ti\{0\}\cup\prt\Gw\ti [0,T)$  and $g\in C(\Gw\ti\BBR_+\times \mathbb R)$. If $u$ is defined in $\Gw\ti\BBR_+$, we denote by $g\!\circ \!u$ the function $(x,t)\mapsto g(x,t,u(x,t))$. We say that $g$ belongs to $\CH$ (resp $\CH_0$) if 
\begin{equation}\label{CH}\BA {ll}
&g(x,t,r)\geq 0\quad\forall (x,t,r)\in \Gw\ti\BBR_+\ti\BBR_+\\[2mm]
&\left(\text{resp. } g\in \CH\text{ and }r\mapsto g(x,t,r)\text{ is nondecreasing}\right.).
\EA\end{equation}

We denote by $\mathfrak M(\Gw)$ the set of Radon measures in $\Gw$, and by $\mathfrak M^b(\Gw)$ (resp.  $\mathfrak M^{b,\gr}(\Gw)) $ the subset of Radon measures such that
$$\myint{\Gw}{}d\abs\gm<\infty\quad\left(\text{resp. }\myint{\Gw}{}\gr d\abs\gm<\infty\right),
$$
where $\gr(x):=\dist(x,\prt\Gw)$. Their positive cones are respectively $\mathfrak M_+(\Gw)$, $\mathfrak M_+^{b}(\Gw)$ and $\mathfrak M_+^{b,\gr}(\Gw)$.

\begin{De} Let $\CS$ be a relatively closed subset of $\Gw$ and $\gm$ a Radon measure on $\CR:=\Gw\setminus\CS$. We say that a nonnegative function $u\in C(Q^\Gw_T)$ admits the couple $(\CS,\gm)$ for initial trace if

\begin{equation*}
\displaystyle \lim_{t\rightarrow 0} \int_{\cal R} u(x,t)\gz(x)dx=\int_{\cal R} \gz d\mu\qquad\forall \gz \in C_c(\cal R),
\end{equation*}
and 
\begin{equation*}
\displaystyle \lim_{t\rightarrow 0} \int_{U} u(x,t)dx=\infty\qquad \forall U \subset \Omega, U \mbox{ open, } U\cap \cal S\neq \emptyset
\end{equation*}
The set $\cal S$ is the set of singular initial points of $u$ and its complement $\cal R$ the set of regular initial points. We write $tr_{\Omega}(u)=(\cal S,\mu)$.
\end{De}

Let $\tilde\gm$ be the extension of $\gm$ as a locally bounded Borel measure. To the couple $(\CS,\tilde\gm)$  we can associate a unique outer regular Borel measure $\gn$ defined by
\begin{equation}\label{borel}
\gn(E)=\left\{\BA {ll}\infty\qquad &\forall E \subset \Omega: E \mbox{ Borel, } E\cap \cal S\neq \emptyset\\
\tilde \gm (E) &\forall E \subset \Omega: E \mbox{ Borel, } E\subset \cal \CR.
\EA\right.
\end{equation}
\blemma{lema}
Assume $\Gw$ is a bounded open domain with a $C^2$ boundary, $T>0$, $g\in \CH$, and let $u\in C(\overline\Gw\ti (0,T])$ be a positive solution of 
\begin{equation}\label{equ}
\prt_tu-\Gd u+g\!\circ \!u=0 \qquad\text{in }Q^\Gw_T.
\end{equation}
If $g\!\circ \!u \in L_\gr^1 (Q^\Gw_T)$, then $u\in L^\infty\left( 0,T, L^1_{\gr}(\Omega)\right)$ and there exists  
$\mu \in \frak M^\gr_+(\Omega)$ such that
\begin{equation}\label{equ}
\lim_{t\to 0}\myint{\Gw}{}u(x,t)\gz (x) dx=\myint{\Gw}{}\gz d\gm\qquad\forall\gz\in C_c(\Gw).
\end{equation}
\es
\Proof 
If $\phi_1 > 0$ is the first eigenfunction of $-\Delta$ in $W^{1,2}_0\left(\Omega\right)$ and $\lambda_1$ is the corresponding eigenvalue, we have
$$\displaystyle \frac{d}{dt}\int_{\Omega} u \phi_1dx+ \lambda_1\int_{\Omega}u \phi_1dx+\int_{\Omega} g\!\circ \!u\, \phi_1 dx+ \int_{\partial \Omega} u \frac{\partial \phi_1}{\partial \nu}dS=0,$$ 
where $\nu$ is the normal vector. Set $X=\myint{\Omega}{}u\phi_1dx$, then by Hopf Lemma,
$$X'+\lambda_1 X + \int_{\Omega} g\!\circ \!u\, \phi_1 dx\geq 0$$
 which yields to
$$\displaystyle \frac{d}{dt} \left( e^{\lambda_1 t }X\right) + e^{\lambda_1 t } \int_{\Omega} g\!\circ \!u\, \phi_1 dx\geq 0.$$
For $s \in (t,T)$
$$\displaystyle \frac{d}{dt}  \left( e^{\lambda_1 t }X - \int_t^T  e^{\lambda_1 s } \int_{\Omega}g\!\circ \!u\, \phi_1 dxds\right)\geq 0,$$ 
which means that the mapping
$$t\rightarrow  e^{\lambda_1 t }X - \int_t^T  e^{\lambda_1 s } \int_{\Omega}g\!\circ \!u\, \phi_1 dxds$$ is nondecreasing.
Therefore
$$e^{\lambda_1 t }X - \int_t^T  e^{\lambda_1 s } \int_{\Omega}g\!\circ \!u\, \phi_1 dxds \leq e^{\lambda_1 T }X.$$
and finally
$$X\leq e^{\lambda_1 (T-t) }X+ e^{-\lambda_1 t } \int_t^T  e^{\lambda_1 s } \int_{\Omega}g\!\circ \!u\, \phi_1 dxds.$$
Since $\gr^{-1}\gf_1$ is positively bounded from above and from below, $u\in L^\infty\left( 0,T, L^1_{\gr}(\Omega)\right)$ and there exists a sequence $\{t_n\}$ decreasing to $0$ and a measure $\gm\in\frak M^\gr_+(\Omega)$ such that
$$\lim_{t_n\to 0}\myint{\Gw}{}u(x,t_n)\gz dx=\myint{\Gw}{}\gz d\gm\qquad\forall\gz\in C_c(\Gw).
$$
If $\gz\in C^2_c(\Gw)$ there holds
$$\myint{\Gw}{}u(x,t_n)\gz dx=\myint{t_n}{T}\myint{\Gw}{}\left(g\!\circ \!u\,\gz-u\Gd\gz\right)dxdt+\myint{\Gw}{}u(x,T)\gz dx,
$$
thus 
$$\myint{\Gw}{}\gz d\gm=\dint_{Q^\Gw_T}\left(g\!\circ \!u\,\gz-u\Gd\gz\right)dxdt+\myint{\Gw}{}u(x,T)\gz dx.
$$
This implies that $\gm$ is uniquely determined and $u(.,t)$ converges to $\gm$ in the weak sense of measures.\qeda

\begin{cor}
Assume $\Gw$ is an open domain, $g\in\CH$ and $u \in C^{2,1}(Q^\Gw_T)$ is a positive solution of (\ref{equ}). Suppose that for any $z\in \Omega$ there exists an open neighborhood $U\subset\Gw$ such that
$$\displaystyle\int_0^T\int_{U} g\!\circ \!u dxdt < \infty.$$
Then $u(x,t)\in L^\infty \left( 0,T,L^1_{loc}(\Omega)\right)$ and there exists a positive Radon measure $\mu$ on $U$ such that 
\begin{equation*}
\displaystyle \lim_{t\rightarrow 0} \int_{\cal R} u(x,t)f(x)dx=\int_{\cal R} f d\mu\quad \forall f \in C_c(\cal R).
\end{equation*}
\end{cor}
\Proof We apply the previous lemma in replacing $U$ by a ball $B_\ge(z)$ and conclude by a partition of unity.\qeda\medskip

The following class of nonlinearity has been introduced by Marcus and V\'eron \cite{MV3}  in order to study the boundary trace of solutions of elliptic equations.

\begin{De}
A function $g\in \CH$ is a coercive nonlinearity in $Q^\Gw_T$ if, for every subdomain $\Omega'$ of $\Omega$ and every $\ge\in (0,T)$, the set of positive solutions of (\ref{0.0}) in $Q^{\Gw'}_{\ge,T}:=\Omega' \times (\ge,T)$ is uniformly bounded in compact subsets of $Q^{\Gw'}_{\ge,T}$.
\end{De}
\begin{De}
Let $z\in \Omega$. We say that equation (\ref{0.0}) possesses a strong barrier at $z$ if there exists a number $r_0 \in (0,\gr(z))$ such that, for every $r \in (0,r_0)$, there exists a positive supersolution $w=w_{r,z}$ of (\ref{0.0}) in $B_r(z)\times (0,T)$ such that 
\begin{equation}\label{barr}
w\in C\left(B_r(z)\times [0,T)\right), \displaystyle \lim _{\left|x-z\right|\rightarrow r} w(x,t)=\infty \mbox { locally uniformly if } t\in (0,T).
\end{equation}
\end{De}

\begin{lema} Assume $g\in\CH$  is a coercive nonlinearity in $Q^\Gw_T$, then the set of solutions of (\ref{0.0}) in $Q^\Gw_T$ is uniformly bounded from above in every compact subset of $Q^\Gw_T$. Furthermore, if $g\in \CH_0$, $A\subset\Omega$ is open and (\ref{0.0}) possesses a strong barrier at every point of $z\in A$, then the set of solutions u of (\ref{0.0}) such that $u \in C\left( A \times [0,T)\right)$ and $u(x,0)=0$ on $A$ is uniformly bounded from above in every compact subset of $A \times [0,T)$.
\end{lema}
\Proof Let $K$ be a compact subset of $Q^\Omega_T$ and let $\Omega'$ be a smooth, bounded domain of $\Gw$ and $\ge>0$ such that $K\subset Q^{\Gw'}_{\ge,T}$ Let $U=U_{Q^{\Gw'}_\infty}$ be the minimal large solution of (\ref{0.0}) in $Q^{\Gw'}_T$, i.e. the limit, when $k\to\infty$, of solutions with Cauchy-Dirichlet data $k$ on $ \prt_\ell Q^{\Gw'}_{\ge,T}:=\Gw'\ti\{\ge\}\cup\prt\Gw'\ti[\ge,T)$. By the maximum principle, if $u\in C(Q^\Omega_T) $ is a solution of (\ref{0.0}), then $u \leq U$ in $\Omega'$.\\
For the second statement,  let $K$ be a compact subset of $A$. For any $z\in K$ there exists $r_z>0$ such that for any $r\in (0,r_z)$ there exists a positive supersolution of (\ref{0.0}) in $Q^{B_r(z)}_T$ which satisfies (\ref{barr}). Since $K$ is compact, there exist
$z_1,...,z_p$ such that $K\subset \cup_{j=1}^pB_{r_{z_j}/2}(z_j)$. For any $j\in \{1,...,p\}$ we denote by $w_j$ the supersolution in $Q^{B_{2r_{z_j}/3}(z_j)}_T$. By comparison principle, there holds
\begin{equation}\label{barr1}
u(x,t)\leq \sup\{w_j(x,t):(x,t)\in B_{r_{z_j}/2}(z_j)\times (0,T)\}:=M_j,
\end{equation}
for $(x,t)\in Q^{B_{r_{z_j}/2}(z_j)}_T$. Therefore $u\leq M=\max_{j=1,...,p}M_j$ in $K\ti (0,T)$.
\qeda
\begin{lema}
Let $g\in\CH$ and $u\in C^{2,1}(Q^\Gw_T)$ be a positive solution of (\ref{0.0}) and suppose $z \in \Gw$ is such that
\begin{equation}\label{sing}
\int_0^T\!\!\int_{B_\ge(z)\cap\Gw}\!\!g\!\circ \!udxdt=\infty\qquad\forall\ge>0.
\end{equation}
Suppose that at least one of the following sets (i) or (ii) of conditions holds:\smallskip

\noindent (i) There exists an open neighborhood $U'\subset\Gw$ of $z$ such that $u\in L^1(U'\times (0,T))$.\smallskip

\noindent (ii) The following hold:

\noindent 1- $g\in \CH_0$,

\noindent 2- (\ref{0.0}) possesses a strong barrier at $z$. \smallskip

\noindent Then,
\begin{equation}\label{sing1}\displaystyle \lim_{t\rightarrow 0} \int_{B_\ge(z)\cap\Gw} u(x,t)dx=\infty\qquad\forall\ge>0.
\end{equation}
\end{lema}
\Proof
Assume that $\Omega$ is bounded. First consider the case when condition (i). holds. Let $\phi \in C^{2,1} \left(U'\times [0,T)\right)$ with compact support in $U'\times [0,T)$ and such that $\phi(x,0)=1$ in a neighborhood of $z$. Then
\begin{equation}\label{sing2}
\displaystyle \int_t^T\!\!\int _{U'}\left(u(-\phi_t-\Delta \phi)+ g\!\circ \!u\,\phi\right)dxdt=\int _ {U'}u(t)\phi dx-\int_{U'}u(T)\phi dx.
\end{equation}
By assumption $\displaystyle \int_t^T\!\!\int _{U'}u(\phi_t+\Delta \phi)dxdt$
 is bounded. We let $t$ tend to 0, the result follows from (\ref{sing}).\smallskip
 
Next we assume that condition (ii) holds, $u\notin  L^1(U'\times (0,T))$ for any neighborhood  $U'$ of $z$ and that the conclusion is not valid. Thus there exist $r^*>0$, such that $\overline B_{r^*}(z)\subset U_z$ and a sequence $\{t_n\}$ decreasing to $0$ such that 
$$ \int_{B_{r^*}(z)} u(x,t_n)dx\leq M$$
for some $M>0$. Furthermore $g$ is coercive in $B_{r^*}(z)\times (0,T)$. Let $\left\lbrace h_{n,k}\right\rbrace \subset C^\infty(Q^\Gw_T)$ an increasing sequence with respect to $k$ and $n$ of nonnegative functions such that $h_{n,k}=0$ on $ B_{r^*}(z) \times\left\lbrace 0\right\rbrace $, $0\leq h_{n,k}\leq k$ and $h_{n,k}=k $ on $(t_n,T)\times \partial B_{r^*}(z)$. Let $w_{h_{n,k}}$ be the solution of  (\ref{0.0}) in $B_{r^*}(z) \times (0,T)$ such that $w_{h_{n,k}}=h_{n,k}$ on $\prt_\ell Q^{B_{r^*}(z)}_{T}$. By the maximum principle and condition (ii)-1, the sequence $\{w_{h_{n,k}}\}$ is monotone increasing  with respect to $k$ and $n$. Condition (ii)-2 implies that, for every $r<{r^*}$ and $\beta<T$, the sequence is bounded in $\overline B_r(z)\times [0,\beta]$, and since $u$ is locally bounded in $Q^\Gw_T$ there exists $k=k(n)$ such that $k\geq u$ on $(t_n,T)\times \partial B_{r}(z)$ and $k(n)\to\infty$ when $n\to\infty$. Then $w =\lim _{n \rightarrow \infty } w_{h_{n,k}}$ is a solution of (\ref{0.0}) which blows up on $\partial B_{r^*}(z) \times (0,T)$ and vanishes on $B_{r^*}(z) \times\{0\}.$
Let $v_{n}$ be the solution of the heat equation in $ B_{r^*}(z)\times (t_n,T)$ such that $v_{n}(.,t_n)=u(.,t_n)$ 
in $B_{r^*}(z)$ and $v_{n}=0$ on $\prt B_{r^*}(z)\ti(t_n,T)$. Then $w_{h_{n,k(n)}}+v_{n}$ is a supersolution of (\ref{0.0}) in $B_{r^*}(z) \times (t_n,T)$ which dominates $u$ on $\prt_\ell Q^{B_{{r^*}}(z)}_{t_n,T}$. By the maximum principle, $$u \leq w_{h_{n,k(n)}}+v_{n}\qquad\text{in }Q^{B_{{r^*}}(z)}_{t_n,T}.$$
And we have in particular
$$\myint{B_r(z)}{}u(x,t)dx\leq \myint{B_r(z)}{}\left(w_{h_{n,k(n)}}+v_{n}\right)(x,t)dx\leq M+\myint{B_r(z)}{}w(x,t)dx\qquad\forall 
t\in (t_n,T).
$$ 
Since it holds for any $n$, it implies $u\in L^1(Q^{B_r(z)}_T)$, which leads to a contradiction. \qeda\medskip

\noindent {\it Example 1}. If $g(x,t,r)=b(x,t)h(r)$ where $b$ is a Borel function defined in $Q^\Gw_T$ which satisfies $\inf\,{\rm \!ess}\,\{b(x)x\in K\}=b_K>0$, $h$ is continuous, nondecreasing and $h(0)\geq 0$, then 
\begin{equation}\label{0}
u_t-\Gd u+h(u)=0
\end{equation}
possesses a strong barrier at any $z\in\Gw$ if and only if $h$ satisfies the Keller-Osserman condition, that is there exists some  $a\geq 0$ such that
\begin{equation}\label{KO}
\myint{a}{\infty}\frac{ds}{\sqrt {H(s)}}<\infty\text{ where }H(s)=\myint{0}{s}h(\gt)d\gt.
\end{equation}
The supersolution can be chosen to be the maximal solution $\gf_r$ of the elliptic equation
\begin{equation}\label{KO1}
-\Gd \gf+b_{\overline B_r(z)}h(\gf)=0\qquad\text{in }B_r(z).
\end{equation}
If we assume moreover that $h$ is super-additive, i.e. $h(a+b)\geq h(a)+h(b)$ for all $a,b\geq 0$, then there holds
\begin{equation}\label{KO}
\myint{a}{\infty}\frac{ds}{h(s)}<\infty,
\end{equation}
and any solution $u$ of  (\ref{0}) is dominated in $Q^{B_r(z)}_T$ by $\gf_r(x)+\psi (t)$ where $\psi$ is defined by inversion from 
$$\myint{\psi (t)}{\infty}\frac{ds}{h(s)}=b_{\overline B_r(z)}t\quad\forall t>0.
$$

\noindent {\it Example 2}. If $g(x,t,r)=a(x)b(t)h(r)$ where $a\in C(\Gw)$, $b\in C((0,T))$, $a,b>0$, then $g$ is a coercive nonlinearity if $h$ is super-additive and satisfies the Keller-Osserman condition. This is not sufficient for the existence of a barrier as it is shown in \cite{MV2} with $h(r)=r^q$ ($q>1$) $a\equiv 1$ and $b(t)=e^{-\frac{1}{t}}$.

\bprop{stab} Let $g\in\CH_0$ such that at any $z\in \Gw$ there exists a strong barrier. We assume also 
\bel{Ap0}
g(x,t,a)+g(x,t,b)\leq g(x,t,a+b)\qquad\forall (x,t,a,b)\in Q_T^\Gw\ti\BBR_+\ti\BBR_+.
\ee
Let $\{u_n\}$ be a sequence of positive solutions of (\ref{0.0}) which converges to $u$ locally uniformly in $Q_T^\Gw$. Denote by  $tr_\Gw(u_n)=(\CS_n,\gm_n)$ and $tr_\Gw(u)=(\CS,\gm)$ their respective initial trace. If $\CA\subset\cap_n\CR_n$ is open and if $\gm_n(\CA_n)$ remains bounded independently of $n\in \BBN$, then $\CA\subset\CR:=\Gw\setminus\CS$.
\es
\Proof Let $z\in \CA$ and $\tilde r\in(0,\gr(z))$ such that for any $r\in (0,\tilde r]$ there exists a positive supersolution $w_{r,z}$ satisfying (\ref{barr}) and $\overline B_{\tilde r}(z)\subset\CA$. For any $n\in\BBN$ and $\gt\in (0,T)$, we denote by $u_{\gt,\chi_{_{B_{\tilde r}(z)}}\gm_n}$ the solution of 
\bel{Ap1}
\BA{ll}
\prt_tu-\Gd u+g\circ u=0\qquad&\text{in } B_{\tilde r}(z)\ti (\gt,T)\\
\phantom{\prt_tu-\Gd u+}
u(.,\gt)=\chi_{_{B_{\tilde r}(z)}}u_n(.,\gt)\qquad&\text{in } B_{\tilde r}(z)\\
\phantom{\prt_t-\Gd u+g\circ u}
u=0\qquad&\text{in } \prt B_{\tilde r}(z)\ti (\gt,T).
\EA\ee
By the maximum principle $u_{\gt,\chi_{_{B_{\tilde r}(z)}}\gm_n}\leq u_n$ in $B_{\tilde r}(z)\ti (\gt,T)$, and 
$g\circ u_{\gt,\chi_{_{B_{\tilde r}(z)}}\gm_n}\leq g\circ u_n $ 
Furthermore, if $\gz\in C^{1,1;1}(Q_T^{\overline B_{\tilde r}(z)})$ vanishes on $\prt B_{\tilde r}(z)\ti [0,T)$ and for $t=T$, there holds
\bel{Ap2}
\dint_{\!\!B_{\tilde r}(z)\ti (\gt,T)}\left(-u_{\gt,\chi_{_{B_{\tilde r}(z)}}\gm_n}(\prt_t\gz+\Gd\gz)+\gz g\circ u_{\gt,\chi_{_{B_{\tilde r}(z)}}\gm_n}\right)dx dt=\myint{B_{\tilde r}(z)}{}\!\!u_n(x,\gt)\gz(x,\gt)dx.
\ee
Since $u_{\gt,\chi_{_{B_{\tilde r}(z)}}\gm_n}$ and $g\circ u_{\gt,\chi_{_{B_{\tilde r}(z)}}\gm_n}$ are bounded independently of $\gt$, standard regularity theory for parabolic equations implies that they converge a.e. in 
$B_{\tilde r}(z)\ti (0,T)$ when $\gt\to 0$ to $u_{\chi_{_{B_{\tilde r}(z)}}\gm_n}$ and $g\circ u_{\chi_{_{B_{\tilde r}(z)}}\gm_n}$. Furthermore
$$\lim_{\gt\to 0}\myint{B_{\tilde r}(z)}{}u_n(x,\gt)\gz(x,\gt)dx=\myint{B_{\tilde r}(z)}{}\gz(x,0)d\gm_n(x).
$$
Using the dominated convergence theorem, it follows from (\ref{Ap1}) that 
\bel{Ap3}
\dint_{\!\!Q_T^{B_{\tilde r}(z)}}\left(-u_{\chi_{_{B_{\tilde r}(z)}}\gm_n}(\prt_t\gz+\Gd\gz)+\gz g\circ u_{\chi_{_{B_{\tilde r}(z)}}\gm_n}\right)dx dt=\myint{B_{\tilde r}(z)}{}\gz(x,0)d\gm_n(x),
\ee
and $u_{\chi_{_{B_{\tilde r}(z)}}\gm_n}$ is the (unique) solution of 
\bel{Ap4}
\BA{ll}
\prt_tu-\Gd u+g\circ u=0\qquad&\text{in } Q_T^{B_{\tilde r}(z)}\\
\phantom{\prt_tu-\Gd u+}
u(.,0)=\chi_{_{B_{\tilde r}(z)}}\gm_n\qquad&\text{in } B_{\tilde r}(z)\\
\phantom{\prt_t-\Gd u+g\circ u}
u=0\qquad&\text{in } \prt B_{\tilde r}(z)\ti (0,T).
\EA\ee
Furthermore, if $\eta$ is the solution of the backward problem
\bel{Ap5}
\BA{ll}
\prt_t\eta+\Gd \eta=-1\qquad&\text{in } Q_T^{B_{\tilde r}(z)}\\
\phantom{\prt_tu}
\!\eta(.,T)=0\qquad&\text{in } B_{\tilde r}(z)\\
\phantom{\prt_t+\Gd \eta}
\eta=0\qquad&\text{in } \prt B_{\tilde r}(z)\ti (0,T),
\EA\ee
there holds
\bel{Ap6}
\dint_{\!\!Q_T^{B_{\tilde r}(z)}}\left(u_{\chi_{_{B_{\tilde r}(z)}}\gm_n}+\eta g\circ u_{\chi_{_{B_{\tilde r}(z)}}\gm_n}\right)dx dt=\myint{B_{\tilde r}(z)}{}\eta(x,0)d\gm_n(x)\leq M,
\ee
for some $M>0$ independent of $n$. Next we set  $Z_{\gt,n}:=u_{\gt\chi_{_{B_{\tilde r}(z)}}\gm_n}+w_{\tilde r,z}$. It is a supersolution of (\ref{0}) in $ (\gt,T)\ti B_{\tilde r}(z)$ which is infinite on $\prt B_{\tilde r}(z)\ti [\gt,T)$ and dominates $u_n$ in $B_{\tilde r}(z)$ at $t=\gt$. Thus $Z_{\gt,n}\geq u_n$ in $ (\gt,T)\ti B_{\tilde r}(z)$. Letting $\gt\to 0$ we finally obtain
\bel{Ap7}
u_{\chi_{_{B_{\tilde r}(z)}}\gm_n}(x,t)\leq u_n(x,t)\leq u_{\chi_{_{B_{\tilde r}(z)}}\gm_n}(x,t)+w_{\tilde r,z}(x,t)\qquad\forall (x,t)\in Q_T^{B_{\tilde r}(z)}.
\ee
For any $r<\tilde r$ and $T'<T$, there exists $\gd, \gs>0$ such that $\eta(x,t)\geq \gd$ and 
$w_{\tilde r,z}(x,t)\leq \gs$ for all $(x,t)\in Q_{T'}^{B_{r}(z)}$. It follows from (\ref{Ap6}), (\ref{Ap7}) and 
Fatou's lemma  that $u$ and $g\circ u$ are integrable in $L^1(Q_{T'}^{B_{r}(z)})$. By \rlemma{lema}, $B_{r'}(z)\subset \CR$. Since it holds for any $z\in\CA$, the result is proved.\qeda\medskip


\mysection{Construction of a barrier}
\def\theequation{3.\arabic{equation}}\makeatother
In the next results we construct the barrier function

\blemma{ss1lem}Assume $\ga>-1$ and $q>1$, then there exists a unique positive function $W_{\ga}\in C^{2}([0,\infty))$ satisfying
\begin{equation}\label{ss2}\BA {l}
W''+\myfrac{r}{2}W'+\myfrac{1+\ga}{q-1}W-W^q=0\quad\text{in }(0,\infty)\\
\phantom{-W''+.\myfrac{1+\ga}{q-1}.}
\lim_{r\to 0}W(r)=\infty\\
\phantom{\myfrac{1+\ga}{q-1}.U-U^q}
\displaystyle\lim_{{r\to\infty}}r^{\frac{2}{q-1}}W(r)=0.
\EA\end{equation}
Furthermore $W_\ga$ is decreasing and 
\begin{equation}\label{ss3}\BA {l}
W_{\ga}(r)=Cr^{\frac{2(1+\ga)}{q-1}-1}e^{-\frac{r^2}{4}}(1+\circ (1))\quad\text{as }\;r\to\infty.
\EA\end{equation}
\es
\Proof Consider the functional 
\begin{equation}\label{ss3}\BA {l}
J( \gf):=\myfrac{1}{2}\myint{0}{\infty}\left( \gf'^2-\myfrac{1+\ga}{q-1} \gf^2+\myfrac{2}{q+1}\abs  \gf^{q+1}\right) e^{\frac{r^2}{4}}dr
\EA\end{equation}
defined over the convex set
$$H_{k}:=\{\gf\in W^{1}_{2}(0,\infty;e^{\frac{r^2}{4}}dr)\cap L^{q+1}(0,\infty;e^{\frac{r^2}{4}}dr): \gf(0)=k\}.
$$
Note that if $\gf\in H_{k}$, 
$$\BA {l}e^{\frac{r^2}{4}}\gf^2(r)=\myint{r}{\infty}(e^{\frac{s^2}{4}}\gf^2(s))'ds\\
\phantom{e^{\frac{r^2}{4}}\gf^2(r)}
=2\myint{r}{\infty}e^{\frac{s^2}{4}}\gf\gf' (s)ds+\myfrac{1}{2}\myint{r}{\infty}se^{\frac{s^2}{4}}\gf^2(s)ds.
\EA$$
In this set $J$ admits a positive minimizer $w_{k}$ which is the unique solution of  
\begin{equation}\label{ssx}\BA {l}
w''+\myfrac{r}{2}w+\myfrac{1+\ga}{q-1}w'-w^q=0\quad\text{in }(0,\infty)\\
\phantom{w''+\myfrac{r}{2}w'\frac{1+\ga}{q-1}w-w^q,}
w(0)=k.
\EA\end{equation}
Furthermore, $w_{k}=\lim_{n\to\infty}w_{{k,n}}$ where $w_{k,n}$ is the unique positive solution of 
\begin{equation}\label{ssxx}\BA {l}
w''+\myfrac{r}{2}w'+\myfrac{1+\ga}{q-1}w-w^q=0\quad\text{in }(0,n)\\
\phantom{w''+\myfrac{r}{2}w'\frac{1+\ga}{q-1}w-w^q,}
w(0)=k\\
\phantom{w''+\myfrac{r}{2}w'\frac{1+\ga}{q-1}w-w^q,}
w(n)=0.
\EA\end{equation}
and, by the maximum principle, $(k,n)\mapsto w_{k,n}$ is increasing. If we consider the linear equation
\begin{equation}\label{ssxy}
z''+\myfrac{r}{2}z'+\myfrac{1+\ga}{q-1}z=0\quad\text{in }(0,\infty),
\end{equation}
it admits two linearly independent positive solutions $z_{1}$ and $z_{2}$ with the following asymptotic behaviour as $r\to\infty$
\begin{equation}\label{ssy1}
z_{1}(r)=r^{-\frac{2(1+\ga)}{q-1}}(1+\circ (1))
\end{equation}
and
\begin{equation}\label{ssy2}
z_{2}(r)=r^{\frac{2(1+\ga)}{q-1}-1}e^{-\frac{r^2}{4}}(1+\circ (1))
\end{equation}
(see \cite[Appendix]{MV1-1}. Since any solution of \ref{ssx}, and \ref{ssxx} as well,  satisfies an a priori estimate of Keller-Osserman type (see \cite{Ve1})
\begin{equation}\label{ssy3}
w(r)\leq Cr^{-\frac{2}{q-1}}\qquad\text{ for } 0<r<1,
\end{equation}
there holds
$$w_{{k,n}}\leq Cz_{2}(r)\qquad\text{ for }  k\geq r\geq 1.
$$
Letting $n$ and $k$ go to infinity successively, it follows that $W_{\ga}=\lim_{k,n\to\infty}w_{k,n}$ exists. It is a positive solution of problem (\ref{ss2}) and it satisfies 
$$W_{\ga}(r)\leq C\left(r^{-\frac{2}{q-1}}+z_{2}(r)\right)\qquad\text{ for } r >0.$$
The singular behaviour at $r=0$ is standard (see e.g. \cite{Ve1}) and yields to
\begin{equation}\label{ssy4}
W_{\ga}(r)=\myfrac{2(q+1)}{(q-1)^2}r^{-\frac{2}{q-1}}(1+\circ (1))\qquad \text{as }r\to 0.
\end{equation}
Thus uniqueness follows by the maximum principle and estimates (\ref{ss3})  is obtained via standard linearization, using the upper estimate at infinity.\qeda\medskip

In the sequel we set

\begin{equation}\label{ssy54}
w_\ga(s,t)=t^{-\frac{1+\ga}{q-1}}W_{\ga}\left(\frac{s}{\sqrt t}\right)\qquad\forall s>0,\,t>0.
\end{equation}

\bprop{est-1} Assume $\ga>-1$ and $q>1$. Then for any $R>0$, there exists $C=C(q,\ga,R)>0$ such that any solution $u$ of (\ref{1}) in $Q^{B_{R}}_{\infty}$ which vanishes on $B_{{R}}\ti\{0\}$satisfies
\begin{equation}\label{ss4}\BA {l}
u(x,t)\leq 2Nt^{-\frac{1+\ga}{q-1}}W_{\ga}\left(\frac{R-\abs{x}}{\sqrt t}\right)\quad\forall (x,t)\in Q^{B_{R}}_{\infty}.
\EA\end{equation}
\es
\Proof For $m>0$, set $S_{m}=\{x=(x_{1},...,x_{N}):\abs {x_{j}}<m,\;\forall j=1,...,N\}$. For $R'<R$
\begin{equation}\label{ss(}\BA {l}
\tilde w_{R'}(x,t)=t^{-\frac{1+\ga}{q-1}}{\displaystyle\sum_{{j=1}}^N}\left(W_{\ga}\left(\frac{R'-x_{j}}{\sqrt t}\right)+W_{\ga}\left(\frac{R'+x_{j}}{\sqrt t}\right)\right)\quad\forall (x,t)\in Q^{S_{R'}}_{\infty}.
\EA\end{equation}
Then $\tilde w_{R'}$ is a supersolution of (\ref{1}) in $Q^{S_{R'}}_{\infty}$ which is infinite on 
$\prt S_{R'}\ti (0,\infty)$, thus $u\leq \tilde w_{R'}$. Letting $R'\to R$ yields to $u\leq \tilde w_{R}$ in $Q^{S_{R}}_{\infty}$. Since the equation is invariant by rotation, for any $x\in B_{R}$, there is a rotation $\CR$ such that $\CR(x)$ has only a positive $x_{1}$- coordinate. Thus 
\begin{equation}\BA {ll}
\label{ssy55}
u(x,t)\leq \tilde w_{R}(\abs x,t)\\[2mm]\phantom{u(x,t)}\leq t^{-\frac{1+\ga}{q-1}}\left(W_{\ga}\left(\frac{R-x_{1}}{\sqrt t}\right)+(2N-1)W_{\ga}\left(\frac{R}{\sqrt t}\right)\right)\\[2mm]\phantom{u(x,t)}
\leq 2Nt^{-\frac{1+\ga}{q-1}}W_{\ga}\left(\frac{R-x_{1}}{\sqrt t}\right),
\EA\end{equation}
which is (\ref{ss4}) since $x_{1}=\abs x$.\qeda\medskip

\bprop{est2} Assume $\ga>-1$, $q>1$ and $R>0$. Then there exists a unique positive solution $w_{B_R}$ of (\ref{1}) in $Q^{B_R}_{\infty}$, continuous in  $B_R\ti [0,\infty)$, which vanishes on $B_R\ti\{0\}$ and satisfies $\lim_{\abs x\to R}w_{B_R}(x,t)=\infty$, locally uniformly in $(0,\infty)$. In particular
\begin{equation}\label{ss4'}\BA {l}
t^{-\frac{1+\ga}{q-1}}W_{\ga}\left(\frac{R-\abs{x}}{\sqrt t}\right)\leq w_{B_R}(x,t)\leq 2Nt^{-\frac{1+\ga}{q-1}}W_{\ga}\left(\frac{R-\abs{x}}{\sqrt t}\right)\quad\forall (x,t)\in Q^{B_{R}}_{\infty}.
\EA\end{equation}
\es
\Proof For $k>0$, let $w^k_{B_{R}}$ be the solution of 
\begin{equation}\label{ss4"}\BA {ll}
\prt_tu-\Gd u+t^\ga u^q=0\qquad&\text {in } Q^{B_{R}}_{\infty}\\\phantom{\prt_t-\Gd u+t^\ga u^q}
u=k\qquad&\text {in } \prt_\ell Q^{B_{R}}_{\infty}\\\phantom{\Gd u+t^\ga u^q}
u(.,0)=0\qquad&\text {in } B_{R}.
\EA\end{equation}
By (\ref{ss4}), $w^k_{B_{R}}(x,t)\leq 2Nt^{-\frac{1+\ga}{q-1}}W_{\ga}\left(\frac{R-\abs{x}}{\sqrt t}\right)$. There there exists 
$w_{B_R}=\lim_{k\to\infty}w^k_{B_{R}}$ and $w_{B_R}$ is a solution of (\ref{1}) in $Q^{B_{R}}_{\infty}$ which vanishes on $B_{{R}}\ti\{0\}$ and is infinite on $\prt B_R\ti (0,\infty)$. Consider the similarity transformation $T_m$ which leaves  equation (\ref{1}) invariant
$$T_m[u](x,t)=m^{\frac{1+\ga}{q-1}}u(\sqrt mx,mt)\qquad\forall m>0,
$$ 
then $T_m[w^k_{B_{R}}]=w_{B_{\frac{R}{\sqrt m}}}^{m^{\frac{1+\ga}{q-1}}k}$ which implies
$$T_m[w_{B_{R}}]=w_{B_{\frac{R}{\sqrt m}}}\qquad\forall m>0.
$$
If $u\in C(B_R\ti [0,\infty))$ is any positive solution of problem
\begin{equation}\label{ss4'"}\BA {ll}
\prt_tu-\Gd u+t^\ga u^q=0\qquad&\text {in } Q^{B_{R}}_{\infty}\\\phantom{t^\ga u^q}
\displaystyle \lim_{\abs x\to R}u(x,t)=\infty\qquad&\text {locally uniformly on } (0,\infty)\\\phantom{\Gd u+t^\ga u^q}
u(.,0)=0\qquad&\text {in } B_{R},
\EA\end{equation}
then for any $m>1$ and $\ge>0$, there exists $\gt_\ge>0$ such that $u(x,t)\leq \ge$ in $B_{\frac{R}{\sqrt m}}$ for $0\leq t\leq \gt_\ge$. Therefore, for any $T>0$, 
$$u\leq \ge+w_{B_{\frac{R}{\sqrt m}}}\qquad \text{on }\prt_\ell Q^{B_{\frac{R}{\sqrt m}}}_T\cup \overline{B_{\frac{R}{\sqrt m}}}\ti\{0\}
$$
Since $\ge+w_{B_{\frac{R}{\sqrt m}}}$ is a supersolution $u\leq \ge+w_{B_{\frac{R}{\sqrt m}}}$ in $Q^{B_{\frac{R}{\sqrt m}}}_T$. Letting $\ge\to 0$, $m\to 1$ and $T\to\infty$ yields to $u\leq w_{B_{R}}$. In the same way, with $0<m<1$, we obtain $u\geq w_{B_{R}}$.\qeda
\medskip

The next estimate is an immediate consequence of \rprop{est-1}.

\bprop{est3} Assume $\ga>-1$, $q>1$ and $K\subset\Gw$ is compact. Let $u$ be any solution $u$ of (\ref{1}) in $Q^{\Gw}_{\infty}$ which vanishes on $\Gw\setminus\{K\}\ti\{0\}$ and on $\prt\Gw\ti [0,\infty)$, then
\begin{equation}\label{ss66}\BA {l}
u(x,t)\leq 2Nt^{-\frac{1+\ga}{q-1}}W_{\ga}\left(\frac{\dist(x,K)}{\sqrt t}\right).
\EA\end{equation}
\es

\medskip

\mysection{Upper estimates}
\def\theequation{4.\arabic{equation}}\makeatother
We start with the following upper estimate already obtained by Shishkov and V\'eron \cite{ShVe} in the case $\alpha \geq 0$.
\bprop{Theo1}
Let $q>1$  and $\alpha >-1$. If $u$ is a solution of (\ref{1}) vanishing on $\prt\Gw\ti [0,T)$, there holds
\begin{equation}\label{es0}
u(x,t)\leq c_{\alpha}t^{-\frac{\alpha+1}{q-1}} \quad\mbox{ for all } (x,t)\in Q^\Gw_T,
\end{equation}
with $c_\alpha=\left(\frac{\alpha+1}{q-1}\right)^{\frac{1}{q-1}}$.
\es
\Proof Let $\phi(t)=c_\alpha t^{-\frac{\alpha+1}{q-1}}$  be the maximal solution of
$$\BA {ll}
\phi' + t^{\alpha}\phi^q=0\\
\phantom{-,-}
\phi(0)=\infty.
\EA
$$
with $c_\alpha=\left(\frac{\alpha+1}{q-1}\right)^{\frac{1}{q-1}}$. \smallskip

\noindent{\it{Case $\alpha \geq 0$.}} For $\tau>0$, we denote by $\Gf_{1,\tau}$ the solution of
\begin{equation}\label{stat}
\begin{array}{lll}
-\Delta  \Gf_{1,\tau}+ \tau^{\alpha} \Gf_{1,\tau}^q=0\quad\mbox{ in }B_1\\
\phantom{-,-}
 \dsp\lim _{|x|\rightarrow 1} \Gf_{1,\tau}(x)=\infty,
\end{array}
\end{equation}
and for $R>0$
$$ \Gf_{R,\tau}(x)=R^{\frac{-2}{q-1}} \Gf_{1,\tau}\left(\frac{x}{R}\right).$$
Note that $\Gf_{R,\tau}(x)$ is the solution of the problem (\ref{stat}) in the ball $B_R$. 
The function $\Gf_{R,\tau}$ tends to 0 uniformly on every compact set of $\mathbb R^N$ when $R\rightarrow \infty$.
Set
$$\tilde v(x,t)=\phi\left(t-\tau\right)+\Gf_{R,\tau}\left(x\right),$$
thus $\tilde v$ is a supersolution of (\ref{1}) in $B_R\times [\tau,T)$ which is infinite  on $\prt B_R\times [\tau,T)\cup B_R\ti\{0\}$. Then $u(x,t)\leq \tilde v(x,t)$.
Letting $R\to\infty$  and $\tau\to 0$, we obtain
$$u(x,t)\leq c_\alpha t^{-\frac{\alpha+1}{q-1}} \quad\mbox{ for all } (x,t)\in Q_T.$$
{\it{Case $-1<\alpha<0.$}}
Let $\tau>0$ and  $\phi_\tau(t)= c_\alpha \left(t^{\alpha+1}-\tau^{\alpha+1}\right)^{-\frac{1}{q-1}}$ be the solution of
$$\BA {ll}
\phi' _\tau+ t^{\alpha}\phi_\tau^q=0\quad\mbox{ on } (\gt,\infty)\\\phantom{-,-}
\phi_\tau(\tau)=\infty
\EA
$$
If  $\Gf_{1,T}$ is the solution of (\ref{stat}) with $\gt=T$, we set
$$\Gf_{R,T}(x)=R^{\frac{-2}{q-1}}\Gf_{1,T}\left(\frac{x}{R}\right).$$ 
Clearly $\Gf_{R,T}$ tends to 0 uniformly on every compact of $\mathbb R^N$ when $R\to \infty$. Set
$$\hat v(x,t)=\phi_\tau (t)+\Gf_{R,T}(x),$$
$\hat v$ is a supersolution of (\ref{1}) in $B_R\times (\tau,T)$, thus $u(x,t)\leq \hat v(x,t)$, as in the first case. 
Letting $R\to\infty$ and $\tau\to 0$, we obtain the desired estimate.\qeda\medskip

Combining \rprop {est-1} and \rprop {Theo1} we obtain,

\bcor{glob} Assume $q>1$, $\ga>-1$ and $K\subset\Gw$ is compact. If $u\in C^{2,1}(Q_T^\Gw)\cap C(\overline Q_T^\Gw\setminus K\ti\{0\})$ is a solution of (\ref{1})  which vanishes on 
$\prt\Gw\ti [0,T)\cup\{(\Gw\setminus K)\ti\{0\}\}$, there holds
\begin{equation}\label{es1}
u(x,t)\leq \min\left\{\dsp 2NW_{\ga}\left({\frac{\dist(x,K)}{\sqrt t}}\right),c_\alpha\right\}t^{-\frac{\alpha+1}{q-1}}
\quad\mbox{ for all } (x,t)\in Q^\Gw_T.
\end{equation}
In the particular case where $K= \{O\}$, (\ref{es1}) yields to
\begin{equation}\label{es2}
u(x,t)\leq \min\left\{\dsp 2NW_{\ga}\left({\frac{\abs x}{\sqrt t}}\right),c_\alpha\right\}t^{-\frac{\alpha+1}{q-1}}\leq \frac{c_1}{\left(\abs x^2+t\right)^{\frac{1+\ga}{q-1}}}
\quad\mbox{ for all } (x,t)\in Q^\Gw_T,
\end{equation}
for some $c_1=c_1(\ga,q)>0$.
\es

\noindent\Remark If $\Gw$ is replaced by $\BBR^N$, the previous estimates (\ref{es0}), (\ref{es1})
 and (\ref{es2}) remain valid. Furthermore, $K$ needs only to be closed.
\mysection{Isolated singularities}
\def\theequation{5.\arabic{equation}}\makeatother
In this section we present the results of classification of isolated singularities of positive solutions of 
(\ref{1}), always in the range $q>1$ and $\ga>-1$. Since some proofs are somewhat similar to the ones of \cite{BF} for the removability of isolated singularities, or \cite{MV1-1}  for the classification of positive  isolated boundary singularities of solutions of 
\begin{equation}\label{sing1}
\prt_tu-\Gd u+\abs u^{q-1}u=0\quad\text{in }Q^\Gw_T,
\end{equation}
we will essentially indicate their main ideas. 
If we look for solution of (\ref{1}) in $Q_T:=Q^{\BBR^N}_T$ under the form
$$u(x,t)=t^\gamma V\left(\frac{x}{\sqrt t}\right)
$$
it is immediate that $\gamma =-\frac{1+\ga}{q-1}$ and $V$ is a solution of 
\begin{equation}\label{sing2}
-\Gd V-\frac{\eta}{2}.\nabla V-\frac{1+\ga}{q-1}V+V^q=0\quad\text{in }\BBR^N.
\end{equation}
It is proved by Escobedo and Kavian \cite{EsKa} that if $\frac{1+\ga}{q-1}>\frac{N}{2}$, or equivalently if
\begin{equation}\label{sing2'}
1<q<q_{c,\ga}:=1+\frac{2(1+\ga)}{N},
\end{equation}
there exists a positive solution of (\ref{sing2}) which minimizes of the functional
\begin{equation}\label{sing3}
\gw\mapsto J(\gw):=\frac{1}{2}\myint{\BBR^N}{}\left(\abs{\nabla \gw}^2-\frac{1+\ga}{q-1}\gw^2+\frac{2}{q+1}\abs\gw^{q+1}\right)e^{\frac{\abs\eta^2}4{}}d\eta,
\end{equation}
over the space $W^{1,2}(\BBR^N;Q d\eta)$ where $Q(\eta)=e^{\frac{\abs\eta^2}4{}}$. The minimizer  $V_\ga$ is unique, radial and satisfies (\ref{sing2}). Furthermore, by adapting the results of \cite{BPT} , it is easy to show that $V_\ga$ is the unique positive $C^2$ function  which satisfies (\ref{sing2}) and 
\begin{equation}\label{sing4-}\BA {ll}
\lim_{\abs\eta\to\infty}\abs\eta^{\frac{2(1+\ga)}{q-1}}V(\eta)=0,
\EA\end{equation}
and that there holds (see \cite[Th 2.1]{MV1-1})
\begin{equation}\label{sing4}
V_\ga(\eta)=C\abs{\eta}^{2\frac{1+\ga}{q-1}-N}e^{-\frac{\abs\eta^2}4{}}(1+\circ (1)) \quad\text{as }\abs\eta\to\infty.
\end{equation}
The function 
\begin{equation}\label{sing5-}
v_\ga(x,t)=t^{-\frac{1+\ga}{q-1}}V_\ga\left(\frac{x}{\sqrt t}\right)
\ee
is a positive solution of (\ref{1}) in $Q_\infty$, continuous in $\overline Q_\infty\setminus\{(0,0)\}$; it vanishes on $\BBR^N\ti\{0\}\setminus\{(0,0)\}$ and satisfies
\begin{equation}\label{sing5}
\lim_{t\to 0}\myint{B_\ge}{}v_\ga(x,t)dx=\infty\qquad\forall\ge>0.
\end{equation}
It is called the {\it very singular solution of (\ref{1})}.\medskip

When $\BBR^N$ is replaced by a a proper open domain $\Gw$ with a compact $C^2$ boundary there exists no self-similar solution to 
(\ref{1}). For any $k>0$ and $a\in\Gw$ there exists a unique solution $u:=u_{k\gd_a}$ to the initial value problem
\bel{init1}\BA {lll}
\prt_tu-\Gd u+t^\ga u^q=0\qquad&\text{in }Q_\infty^\Gw\\
\phantom{\prt_tu-\Gd +t^\ga u^q}
u=0\qquad&\text{on }\prt\Gw\ti [0,\infty)\\
\phantom{-\Gd +t^\ga u^q}
u(.,0)=k\gd_a &\text{in }\Gw.
\EA\ee
(see e.g. \cite{MV2}). The function $u$ belongs to $L^q(Q_T^\Gw;t^\ga dx dt)\cap L^1(Q_T^\Gw)$,  $T>0$ arbitrary, and satisfies
\bel{init2}\BA {lll}
\dint_{\!\!Q_T^\Gw}\left(-u\left(\prt_t\gz+\Gd\gz\right)+t^{\ga}u^q\gz\right) dx dt=k\gz(a)
\EA\ee
for all $\gz\in C^{1,1;1}(\overline Q_T^\Gw)$ which vanishes on $\prt\Gw\ti [0,T]$ and on $\Gw\ti \{T\}$. It is unique in the class of of weak solutions, i.e. the functions belonging to $L^q(Q_T^\Gw;t^\ga dx dt)\cap L^1(Q_T^\Gw)$ and satisfying the above identity. When $k\to\infty$, 
$u_{k\gd_a}\uparrow u_{\infty\gd_a}$, where $u_{\infty\gd_a}:=u_{\infty,a}$ is a solution of (\ref{1}) in $Q_\infty^\Gw$ which vanishes on $\prt\Gw\ti [0,\infty)$ and on $\overline\Gw\ti \{0\}\setminus \{(a,0)\}$ and satisfies (\ref{sing5}). Finally, if $E(x,t)=(4\gp t)^{-\frac{N}{2}}e^{-\frac{\abs x^2}{4t}}$ denotes the heat kernel in $\BBR^N$,
\bel{init3}\BA {lll}
u_{k\gd_a}\sim kE(x-a,t)\qquad\mbox{ when }(x,t)\to (a,0)
\EA\ee
and 
\bel{init4}\BA {lll}
u_{\infty,a}\sim v_\ga(x-a,t)\qquad\mbox{ when }(x,t)\to (a,0).
\EA\ee\smallskip

The following classification of isolated singularities holds.
\bth{class} Assume $\ga>-1$, $1<q<q_{c,\ga}$ and $a\in\Gw$. If $u\in \overline Q^\Gw_T\setminus\{(a,0)\}$ is a positive
solution of (\ref{1}) which vanishes on $\Gw\ti\{0\}\setminus\{(a,0)\}$, then\smallskip

\noindent (i) either there exists $k\geq 0$ such that 
\begin{equation}\label{sing6}
u(x,t)\sim kE(x-a,t)\qquad\mbox{ when }(x,t)\to (a,0),
\end{equation}
and $u$ is a solution of
\begin{equation}\label{sing7}\BA {ll}
\prt_tu-\Gd u+t^\ga u^q=0\qquad&\text{in }Q_T^\Gw\\
\phantom{-----}
u(.,0)=k\gd_{a}\qquad&\text{in }\Gw,
\EA\end{equation}
\smallskip

\noindent (ii) or 
\begin{equation}\label{sing7}
u(x,t)\sim v_\ga(x-a,t)\qquad\mbox{ when }(x,t)\to (a,0),
\end{equation}
and $u$ satisfies (\ref{sing5}).
\es
\Proof The initial trace $tr_\Gw(u)$ is a Borel measure concentrated at $a$ and either it is of the form $(\{\emptyset\},k\gd _a)$ for some $k\geq 0$ or of the form $(\{a\},0)$. In the first case  $t^\ga u^q\in L^1(Q_T^{B_r(a)})$ for any $0<r<\gr(a)$. For $t>0$ we set $m_r(t)=\max\{u(x,\gt):(x,\gt)\in \prt_\ell Q_t^{B_r(a)}\}$, 
$e_r(t)=\max\{E(x,\gt):(x,\gt)\in \prt_\ell Q_t^{B_r(a)}\}$. We denote by $\Gth_r$ the solution of
$$\BA {ll}\prt_\gt\Gth_r-\Gd \Gth_r=0\qquad\text{in }Q_\infty^{B_r(a)}\\
\phantom{\Gd \Gth_r}
\Gth_r(x,0)=0\qquad\text{in }B_r(a)
\\
\phantom{\Gd \Gth_r}
\Gth_r(x,t)=1\qquad\text{in }\prt_\ell Q_\infty^{B_r(a)}
\EA$$
 and
$$\Psi(x,t)=\myint{0}{t}\myint{B_r(a)}{}E(x-y,t-s)s^\ga u^q(y,s)dy ds.
$$
Then
\begin{equation}\label{sing8}
kE(x-a,t)-\Psi(x,t)-e_r(t)\Gth_r(x,t)\leq u(x,t)\leq kE(x-a,t)+m_r(t)\Gth_r(x,t).
\end{equation}
Using the explicit expression of the Cauchy-Dirichlet heat kernel in $B_r$, one can easily check that
$\lim_{t\to 0}\Gth_r(x,t)=0,
$
uniformly on compact subsets of $B_r(a)$. Furthermore
\begin{equation}\label{sing9}\BA {l}\Psi(x,t)\leq 2^{q-1}\myint{0}{t}\myint{B_r(a)}{}E(x-y,t-s)s^\ga \left(k^qE^q(y-a,s)+m^q_r(s)\Gth^q_r(y,s)\right)dy ds\\[4mm]\phantom{\Psi(x,t)}
\leq 2^{q-1}k^q\myint{0}{t}\myint{B_r(a)}{}E(x-y,t-s)s^\ga E^q(y-a,s)dy ds+o(1),
\EA\end{equation}
since the second term in the integral is bounded. Furthermore the first term of the right-hand side of 
(\ref{sing9}) converges to $0$ in $L^1(B_r(a))$ when $t\to 0$. Since $y\mapsto E^q(y-a,s)$ is radially decreasing with respect to $a$, it implies that $\myint{0}{t}\myint{B_r(a)}{}E(x-y,t-s)s^\ga E^q(y-a,s)dy ds$
is maximal at $x=0$ and therefore $\Psi(x,t)\to 0$, uniformly in $B_r(a)$. It follows from (\ref{sing8}) that 
\begin{equation}\label{sing10}
\abs{u(x,t)-kE(x-a,t)}\to 0\qquad\text{when }t\to 0,
\end{equation}
uniformly on compact subsets of $B_r(a)$, for any $r<\gr(a)$.\smallskip

Next we assume that $tr_{\Gw}(u)=(\{a\},0)$, and without loss of generality, we can suppose that $a=0$ and set $B_r=B_r(0)$. Then  
\begin{equation}\label{sing11}
u_{\infty,0}-m_r(t)\Gth_r(x,t)\leq u(x,t)\qquad\forall (x,t)\in Q^{B_r}_\infty,
\end{equation}
where $u_{\infty,0}=\lim_{k\to\infty}u_{k\gd_0}$, and $u_{k\gd_0}$ is the solution of (\ref{1}) in $Q_\infty$ with initial data $k\gd_0$. Although estimate (\ref{sing11}) is proved in \rth{sub}, in next section, its proof does not require any element of the proof of the present theorem. Moreover, for any $0<\ge<r$, 
\begin{equation}\label{sing12}
 u(x,t)\leq V_\ge+m_r(t)\Gth_r(x,t)\qquad\forall (x,t)\in Q^{B_r}_\infty,
\end{equation}
where $V_\ge$ is the limit, when $\ell\to\infty$ of the solutions $V_{\ell,\ge}$ of (\ref{1}) in $Q_\infty$ which has initial data $m\chi_{_{\overline B_\ge}}$. Consider the similarity transformation $T_m[\gf](x,t)=m^{-\frac{1+\ga}{q-1}}\gf(\frac{x}{\sqrt m},\frac{t}{m})$ ($m>0$) which leaves (\ref{1}) and $Q_\infty$ invariant, then, by uniqueness
\begin{equation}\label{sing13}
T_m[u_{k,0}]=u_{m^{\frac{N}{2}-\frac{1+\ga}{q-1}}k,0}
\end{equation}
and
\begin{equation}\label{sing14}
T_m[V_{\ell,\ge}]=V_{m^{\frac{N}{2}-\frac{1+\ga}{q-1}}\ell,\sqrt m\ge}.
\end{equation}
Letting $k,\ell\to\infty$ and $\ge\to 0$, we obtain that, for any $m>0$, 
\begin{equation}\label{sing15}
T_m[u_{\infty,0}]=u_{\infty,0}
\end{equation}
and
\begin{equation}\label{sing16}
T_m[V_0]=V_0.
\end{equation}
Therefore $u_{\infty,0}$ and $V_0$ are positive self-similar solutions of (\ref{1}) in $Q_\infty$ with initial trace $(\{0\},0)$. Then they coincide with the function $v_\ga$ defined in (\ref{sing5-}). The result follows from this equality since (\ref{sing11}), jointly with (\ref{sing12}), implies
\begin{equation}\label{sing17}
v_\ga(x,t)-m_r(t)\Gth_r(x,t)\leq u(x,t)\leq v_\ga(x,t)+m_r(t)\Gth_r(x,t)\qquad\forall (x,t)\in Q^{B_r}_\infty.
\end{equation}
\qeda

\medskip
If $q\geq q_{c,\ga}$ there holds the following result which extends Brezis and Friedman's classical one. 
\bth{remov} Assume $\ga>-1$, $q\geq q_{c,\ga}$ and $a\in\Gw$. If $u\in \overline Q^\Gw_T\setminus\{(a,0)\}$ is a positive solution of (\ref{1}) which vanishes on $\Gw\ti\{0\}\setminus\{(a,0)\}$, then it can be extended as a continuous function $\tilde u$ which vanishes on $\Gw\ti\{0\}$.
\es
\Proof Up to modifying a few parameters the proof 
is similar to Brezis-Friedman's construction. The first step is to prove that $u\in L^q (Q_T^{{B_R}})$ for some $R>0$. This is done by using (\ref{es2}) and the same test function used in \cite{BF}. As a consequence  we obtain that $u$ satisfies 
$$\lim_{t\to 0}\int_{\Gw}u(x,t)\gz(x) dx=0\qquad\forall \gz\in C^{\infty}_0(\Gw).
$$
Finally the extension of $u$ by zero at $t=0$ satisfies the equation in $\Gw\ti[0,T)$.\qeda
\medskip

\noindent\Remark We recall that $E(x,t)=(4\gp t)^{-\frac{N}{2}}e^{-\frac{\abs x^2}{4t}}$. Then if $1\leq r<q_{c,\ga}$, there holds
\bel{int1}
\dint_{Q_T}E^r(x,t) t^\ga dx dt<\infty,
\ee
while if $r\geq q_{c,\ga}$
\bel{int2}
\dint_{Q_T}E^r(x,t) t^\ga dx dt=\infty.
\ee
\mysection{The trace theorem}
\def\theequation{6.\arabic{equation}}\makeatother

In all this section we assume that $\Gw\subset \BBR^N$ is an open domain with a compact $C^2$ boundary, $\ga>-1$ and $1<q<q_{c,\ga}$, and let $u\in C(\overline\Gw\ti (0,T])$ be a positive solution of (\ref{1}) in $Q^\Gw_T$ which vanishes on $\prt\Gw\ti (0,T]$. By Section 2 $u$ possesses an initial trace $tr_\Gw(u)=(\CS,\gm)$ where $\CS$ is a relatively closed subset of $\Gw$ and 
$\gm$ is a Radon measure on $\CR:=\Gw\setminus \CS$.  To this couple we can associate a unique 
outer regular Borel measure $\gn$ defined by
\bel{bo1}
\gn (E)=\left\{\BA {ll}
\gm (E)\qquad&\text{if }E\subset\CR\\
\infty\qquad&\text{if }E\cap\CS\neq\emptyset
\EA\right.
\ee
for any Borel subset $E$ of $\Gw$. Conversely, to any outer Borel measure $\gn$ on $\Gw$ we can associate the regular set $\CR\subset\Gw$ which is the set of points $y\in\Gw$ which possess an open neighborhood $\CO_y$ such that $\gn(\CO_y)<\infty$. Clearly $\CR$ is open and the restriction of $\gn$ to $\CR$ is a positive Radon measure. The set $\CS=\Gw\setminus\CR$ is relatively closed and it is the singular part of $\gn$. It has the property that $\gn(E)=\infty$ for any Borel set $E$ such that $E\cap\CS\neq\emptyset$. We shall denote by $\frak B^{reg}(\Gw)$ the set of outer regular Borel measures in $\Gw$ and by $\frak B^{reg}_c(\Gw)$ the subset of $\frak B^{reg}(\Gw)$ for which $\CS$ is a compact subset of $\Gw$. Thus  $u$ is a solution of the following problem,
\bel{bo1}
\BA {ll}
\prt_tu-\Gd u+t^\ga u^q=0\qquad&\text{in }Q^\Gw_T\\
\phantom{---t^\ga}
u\geq 0,\,u=0\qquad&\text{on }\prt\Gw\ti [0,T)\\
\phantom{---,,t^\ga}
tr_\Gw(u)=\gn\approx (\CS,\gm)\in \frak B^{reg}(\Gw).
\EA\ee

\bdef {CU}We denote by $\CU_{\CS,\gm}(\Gw)$ the set of solutions of problem (\ref{bo1}).
\es

 The first step in the characterization of the singular set is the following delicate lower estimate.

\bth{sub} Let $u\in \CU_{\CS,\gm}(\Gw)$ and $a\in\CS$, then 
\bel{sub1}
u(x,t)\geq u_{\infty,a}(x,t).
\ee
Furthermore, if $\CS$ has a nonempty interior $\CA$, there holds
\bel{sub2}
\lim_{t\to 0}t^{\frac{1+\ga}{q-1}}u(x,t)=c_\ga,
\ee
uniformly on compact subsets of $\CA$.
\es

We first give a proof of (\ref{sub1}) in the case where either $-1<\ga\leq 0$, or $\ga>0$ and $1<q<q_{c,0}$.
\bprop{sub-1} Assume either $-1<\ga\leq 0$ and  $1<q<q_{c,\ga}$, or $\ga>0$ and $1<q<q_{c,0}$, then inequality (\ref{sub1}) holds.
\es
\Proof For any $\ge>0$, there holds
\bel{buv}\lim_{t\to 0}\myint{B_\ge(a)}{}u(x,t) dx=\infty.
\ee
If $k>0$ is fixed, and $\{t_n\}$ is a sequence decreasing to $0$. There exists $t_{n_1}$ such that 
$$\myint{B_{2^{-1}}(a)}{}u(x,t_{n_1}) dx>k,
$$
and there exists $m=m_1(k)$ such that 
$$\myint{B_{2^{-1}}(a)}{}m_1(k)\wedge u(x,t_{n_1}) dx=k,
$$
where $A\wedge B=\inf\{A,B\}$. Assume we have constructed $t_{n_j}<t_{n_{j-1}}$ and $m_j(k)>0$ such that
$$\myint{B_{2^{-j}}(a)}{}m_j(k)\wedge u(x,t_{n_j}) dx=k.
$$
Since (\ref{buv}) holds with $\ge=2^{-j-1}$, there exists $t_{n_{j+1}}<t_{n_j}$ such that 
$$\myint{B_{2^{-j-1}}(a)}{}u(x,t_{n_j}) dx>k.
$$
and thus $m_{j+1}(k)>0$ such that
$$\myint{B_{2^{-j-1}}(a)}{}m_{j+1}(k)\wedge u(x,t_{n_{j+1}}) dx=k.
$$
Next we denote by $u_j$ the solution of 
\bel{bo2}
\BA {ll}
\prt_tu-\Gd u+t^\ga u^q=0\qquad&\text{in }\Gw\ti (t_{n_j},T)\\
\phantom{---t^\ga u\geq 0,\,}
u=0\qquad&\text{on }\prt\Gw\ti [t_{n_j},T)\\
\phantom{---t^\ga}
u(.,t_{n_j})=m_{j+1}(k)\wedge u(.,t_{n_{j+1}})&\text{in } \Gw.
\EA\ee
By the maximum principle $u_j\leq u$ in $\Gw\ti [t_{n_j},T)$, or equivalently 
$u_j(x,t+t_{n_j})\leq u(x,t+t_{n_j})$  in $Q_{T-t_{n_j}}^\Gw$. Clearly
$$u_{j,0}:=u(.,t_{n_j})\rightharpoonup k\gd_a\quad\text{as }j\to\infty
$$
in the weak sense of measures. In order to prove that $u_j$ converges to $u_{k\gd_a}$ we notice that $u_j(x,t+t_{n_j})\leq E\ast u_{j,0}(x,t)$ in $Q_\infty^\Gw$.
If $-1<\ga\leq 0$ and $r\in (1,q_{c,\ga})$ 
$$\BA {ll}\dint_{\!\!Q_T^\Gw}{}u^r_j(x,t+t_{n_j})(t+t_{n_j})^\ga dx dt
\leq \dint_{\!\!Q_T}{}(E\ast u_{j,0})^r(x,t)(t+t_{n_j})^\ga dx dt\\[4mm]
\phantom{\dint_{\!\!Q_T^\Gw}{}u^r_j(x,t+t_{n_j})(t+t_{n_j})^\ga dx dt}
\leq \norm{u_{j,0}}_{L^1}^r\dint_{\!\!Q_T}{}E^r(x,t)t^\ga dx dt,
\EA
$$
 using Young's inequality and (\ref{int1}). If $\ga\geq 0$ and $r\in (1,q_{c,0})$
 $$\BA {ll}\dint_{\!\!Q_T^\Gw}{}u^r_j(x,t+t_{n_j})(t+t_{n_j})^\ga dx dt
\leq 2^\ga T^\ga\dint_{\!\!Q_T}{}(E\ast u_{j,0})^r(x,t) dx dt\\[4mm]
\phantom{\dint_{\!\!Q_T^\Gw}{}u^r_j(x,t+t_{n_j})(t+t_{n_j})^\ga dx dt}
\leq 2^\ga T^\ga\norm{u_{j,0}}_{L^1}^r\dint_{\!\!Q_T}{}E^r(x,t) dx dt.
\EA
$$
 Furthermore, for $s\in (1,q_{c,0})$
$$\BA {ll}\dint_{\!\!Q_T^\Gw}{}u^s_j(x,t+t_{n_j}) dx dt
\leq \dint_{\!\!Q_T}{}(E\ast u_{j,0})^s(x,t) dx dt\\[4mm]
\phantom{\dint_{\!\!Q_T^\Gw}{}u^s_j(x,t+t_{n_j}) dx dt}
\leq \norm{v_{j,0}}^s_{L^1}\dint_{\!\!Q_T}{}E^s(x,t) dx dt.
\EA
$$
Thus the sets of functions $\{u^q_j(.,.+t_{n_j})(.+t_{n_j})^\ga\}$ and $\{u_j(.,.+t_{n_j})\}$ are uniformly integrable in  $L^1(Q_T^\Gw)$. Since $u_j$ satisfies 
\bel{bo4}
\dint_{\!\!Q_T^\Gw}\left(-u_j(x,t+t_{n_j})\left(\prt_t\gz+\Gd\gz\right)+(t+t_{n_j})^{\ga}u^q_j(x,t+t_{n_j})\gz\right) dx dt=\myint{\Gw}{}u_{j,0}\gz dx
\ee
for any test function $\gz\in C^{1,1;1}(\overline Q_T^\Gw)$ vanishing on $\prt\Gw\ti [0,T]\cup \Gw\ti\{T\}$ and converges a.e. in $Q_T^\Gw$ to some $u^*$ when $j\to\infty$, it follows by Vitali's theorem that
\bel{bo4}
\dint_{\!\!Q_T^\Gw}\left(-v\left(\prt_t\gz+\Gd\gz\right)+t^{\ga}v^q\gz\right) dx dt=k\gz(a)
\ee
thus $u^*=u_{k\gd_a}$ by uniqueness, which implies the claim since $u\geq u^*$.\qeda\medskip

When $\ga>0$ and $q_{c,0}\leq q<q_{c,\ga}$, this argument cannot work since the sequence $u_{j,0}$ could concentrate too quickly with respect to $t$ to a Dirac mass and isolated singularities are removable for solutions of 
$$\prt_tu-\Gd u+cu^q=0.
$$
We develop below a proof  valid for any $\ga>-1$ and $1<q<q_{c,\ga}$, which is based upon the parabolic Harnack inequality and shows that such a fast concentration never occurs. 

\blemma{lem-1}
Assume $\ga>-1$ and $1<q<q_{c,\ga}$. Let $\{(x_n,t_n)\}\subset Q_T^\Gw$ be a sequence converging to $(a,0)$ and $\ell>0$. Put $V_n=B_{\ell \sqrt{t_n}}(x_n)$ and suppose that there exist nonnegative functions $h_n\in L^\infty(\BBR^N)$ with support in  $V_n$ such that $0\leq h_n\leq c_1t_n^{-{\frac{N}{2}}}$ and 
\bel{bo5}
h_n\rightharpoonup k\gd_0.
\ee
Then the solutions $u_n$ of 
\bel{bo2}
\BA {ll}
\prt_tu-\Gd u+t^\ga u^q=0\qquad&\text{in }\Gw\ti (t_{n},T)\\
\phantom{---t^\ga u\geq 0,\,}
u=0\qquad&\text{on }\prt\Gw\ti [t_{n},T)\\
\phantom{--,,-t^\ga}
u(.,t_{n})=h_n&\text{in } \Gw,
\EA\ee
satisfy $u_n\to u_{k\gd_a}$ when $n\to\infty$.
\es
\Proof  The estimate $h_n\leq ct_n^{-{\frac{N}{2}}}\chi_{_{V_n}}$ can be written under the form
$$h_n(x)\leq c_2E(x-x_n,t_n)\chi_{_{V_n}}
$$
where $c_2=(4\gp)^{{\frac{N}{2}}}e^{\frac{\ell}{4}}c_1$. 
By the maximum principle
$$u_n(x,t)\leq c_2E(x-x_n,t)\quad\text{ in }\Gw\ti (t_n,\infty).
$$
By (\ref{int1}), (\ref{int2}), the sequences $\{E^q(.-x_n,.)t^\ga\}$ and $\{E(.-x_n,.)\}$ are uniformly integrable in $Q_T^\Gw$, therefore, if we extend $u_n$ by zero in $Q^\Gw_{t_n}$ and denote by 
$\tilde u_n$ the extended function defined in $Q^\Gw_{T}$, we infer that the sequences $\{t^\ga\tilde u^q_n\}$ and $\{\tilde u_n\}$ are uniformly integrable in $Q^\Gw_{T}$. Using standard regularity estimates there exists 
a function $u^*$ defined in $Q^\Gw_{T}$ and a subsequence $u_{n_j}$ such that 
$\tilde u_{n_j}\to u^*$ locally uniformly in $Q^\Gw_{T}$. It follows by uniform integrability and Vitali's convergence theorem that 
$$\tilde u_{n_j}\to u^*\qquad\text{in }L^q(Q^\Gw_{T}; t^\ga dx dt)\cap L^1(Q^\Gw_{T}).
$$
Let $\gz\in C^{1,1;1}(\overline Q^\Gw_{T})$ vanishing on $\prt\Gw\ti [0,T]\cup \Gw\ti\{T\}$, then
$$\dint_{\!\!Q_T^\Gw}\left(-\tilde u_{n_j}(\prt_t\gz+\Gd\gz)+t^\ga \tilde u^q_{n_j}\gz \right)dx dt=\myint{\Gw}{}\gz(.,t_{n_j}) h_n dx.
$$
Using the previous convergence results and the assumption (\ref{bo5}), we derive
$$\dint_{\!\!Q_T^\Gw}\left(-u^*(\prt_t\gz+\Gd\gz)+t^\ga u^{*q}\gz \right)dx dt=k\gz(a).
$$
Thus $u^*=u_{k\gd_a}$ and $\tilde u_{n}\to u_{k\gd_a}$ locally uniformly in $Q^\Gw_{T}$.\qeda\medskip


\blemma{harn} Let $u$ be a positive solution of (\ref{1}) in $Q_T^\Gw$ vanishing on $\prt\Gw\ti [0,T)$. Then for any $\Gw'\subset\overline\Gw\subset\Gw$ there exists a constant $C=C(N,q,\ga,\Gw')>0$ such that 
\bel{h1}
u(y,s)\leq u(x,t)e^{C\left(\frac{\abs{x-y}^2}{t-s}+\frac{t}{s}+1\right)}\qquad\forall (x,t), (y,s)\in Q^{\Gw'}_T
,\, s<t.\ee
\es
\Proof By (\ref{es0}), $V(x,t):=t^\ga u^{q-1}\leq c_\ga^{q-1}t^{-1}$. If we write (\ref{1}) under the form
\bel{h2}
\prt_t u-\Gd u+V(x,t) u=0\quad\text {in } Q^{\Gw}_T
\ee
it follows (\ref{h1}) from parabolic Harnack inequality (see e.g. \cite[Lemma 3.16]{GkVe} although the result is much older).\qeda\medskip

If $G\subset\BBR^N$ is a bounded open subset, we denote by $\gf_G$ is the first eigenfunction
 of $-\Gd$ in $W^{1,2}_0(G)$ normalized by $\sup \gf_G=1$ with corresponding eigenvalue $\gl_G$. 
\blemma {trl}Let $G\subset\BBR^N$ be a bounded open subset with a smooth boundary $\gw\in C^1(Q^{G}_T)$, $\gw\geq 0$,  such that 
\bel{h3}
\myint{0}{T}\norm{\gw(.,t)}_{L^\infty(G)} dt<\infty.
\ee
If $v\in C^{2;1}(\overline G\ti (0,T])$ is a positive solution of 
\bel{h2}
\prt_t v-\Gd v+\gw v=0\quad\text {in } Q^{G}_T,
\ee
then $v\gf_G\in L^1(Q_T^G)$, $\gw v\gf^3_G\in L^1(Q_T^G) $ and there exists $\gm\in\frak M_+(G)$ such that 
\bel{h4}
\lim_{t\to 0}\myint{G}{}v(x,t)\gz (x)dx=\myint{G}{}\gz d\gm\qquad\forall\gz\in C_c{(\overline G}).
\ee
\es
\Proof Set $\gamma (t)=\norm{\gw(.,t)}_{L^\infty(G)}$. Denote $X(t)=\int_{G}{}v(x,t)\gf_G (x)dx$. 
Then, from (\ref{h2}),
 $$X'+\gl_GX+
 \gamma (t)X\geq 0
 $$
 This implies that the function
 $$t\mapsto e^{t\gl_G+\int_0^t\gamma (s) ds}X(t)
 $$
 (which exists thanks to (\ref{h3})) is nondecreasing. Therefore there exists $X(0)=\lim_{t\to 0}X(t)$ and $v\gf_G\in L^1(Q_T^G) $. Furthermore, if we set $Y(t)=\int_{G}{}v(x,t)\gf^3_G (x)dx$
 $$\myfrac{d}{dt}(e^{3t\gl_G}X(t))+e^{3t\gl_1}\int_{\Gw}(\gw\gf^3_G-6\gf_G\abs{\nabla\gf_G}^2)v dx=0.
 $$
Since $\nabla \gf_G$ is bounded and $v\gf_G\in L^1(Q_T^G) $, it implies

 $$e^{3t\gl_G}Y(T)+\int_0^Te^{3t\gl_G}\int_{\Gw}\gw(x,t) v(x,t) \gf^2_G dx dt=Y(0)
 +6\int_0^Te^{3t\gl_G}\int_{\Gw}\gf_G\abs{\nabla\gf_G}^2)v dx dt,
 $$
 which implies that $\gw v\gf^3_G\in L^1(Q_T^G) $. The argument given in the proof of \rlemma{lema} shows that $v$ admits an initial trace which belongs to $\frak M_+(G)$.\qeda\medskip
 
\noindent {\it Proof of \rth{sub}}. We define the parabolic distance in $\BBR^N\ti\BBR$ by
$$d_P((x,t),(y,s))=\sqrt{\abs{x-y}^2+\abs{t-s}}.
$$
{\it Step 1. }We first prove that if $u$ satisfies
\bel{h5}
\limsup_{d_P((x,t))(a,0))\to 0}t^{\frac{N}{2}}u(x,t)<\infty,
\ee
then $a\in\CR(u)$. If (\ref{h5}) holds there exists $\ge,c>0$ such that 
$$u(x,t)\leq ct^{-\frac{N}{2}}\qquad\forall (x,t)\text{ s.t. } \sqrt{\abs{x-a}^2+t}\leq \ge.
$$
If we set $\gw(x,t)=t^\ga u^{q-1}(x,t)$, then 
$$\gw(x,t)\leq c^{q-1}t^{-(\frac{N(q-1)}{2}-\ga)}\qquad\forall (x,t)\in B_{\frac{\ge}{\sqrt 2}}(a)\ti 
(0,\frac{\ge^2}{2}].
$$
Since $q<q_{c,\ga}$, then $\frac{N(q-1)}{2})-\ga<1$; thus the assumptions of \rlemma{trl} are fulfilled and there exists a positive Radon measure $\gm$ in $B_{\frac{\ge}{\sqrt 2}}(a)$ such that 
$$\myint{\Gw}{}u(x,t)\gz(x) dx\to \myint{\Gw}{}\gz d\gm\,\text{ when }t\to 0,\quad\forall \gz\in C^\infty_c(B_{\frac{\ge}{\sqrt 2}}(a)).
$$
Furthermore $t^\ga u^q\in L^1(B_{\frac{\ge}{\sqrt 3}}(a)\ti (0,T))$, which is the claim.  \smallskip

\noindent {\it Step 2. }Since $a\in\CS(u)$, there holds
\bel{h6}
\limsup_{d_P((x,t))(a,0))\to 0}t^{\frac{N}{2}}u(x,t)=\infty.
\ee
Then there exists a sequence $\{(x_n,s_n)\}$ converging to $(a,0)$ such that 
\bel{h7}
u(x_n,t_n)\geq nt_n^{-\frac{N}{2}}.
\ee
We apply \rlemma{harn} with $s=s_n$, $t=2s_n:=t_n$, $y=x_n$ $\abs {x-x_n}\leq\sqrt{ s_n}$.
Then 
$$u(x,t_n)\geq Cnt_n^{-\frac{N}{2}}\qquad\forall x\in V_n:= B_{\frac{\sqrt{t_n}}{2}}(x_n).
$$
This implies 
$$\myint{V_n}{}u(x,t_n) dx\geq C_Nn.
$$
For $k<n$ fixed, we denote by $v:=v_{n,k}$ the solution of 
\bel{h8}\BA {ll}
\prt_tv-\Gd v+t^\ga v^q=0\qquad&\text{in }\Gw\ti (t_n,\infty)\\
\phantom{\prt_t-\Gd v+t^\ga w^q}
v=0\qquad&\text{in }\prt\Gw\ti (t_n,\infty)\\
\phantom{\Gd v,,t^\ga w^q}
v(.,t_n)=Ckt_n^{-\frac{N}{2}}\chi_{V_n}\qquad&\text{in }\Gw.
\EA\ee
By the maximum principle $u\geq v_{n,k}$ in $\Gw\ti (t_n,\infty)$. Furthermore
$$\myint{V_n}{}v(.,t_n) dx=C_Nk.
$$
Thus $v(.,t_n) \rightharpoonup C_Nk\gd_a$ in the weak sense of measures. It follows from \rlemma{lem-1} that $v_{n,k}\to u_{C_Nk\gd_a}$ locally uniformly in $Q_T^\Gw$. Therefore 
$u\geq u_{C_Nk\gd_a}$ in $Q_T^\Gw$. Since $k$ is arbitrary, we obtain (\ref{sub1}).\smallskip

\noindent{\it Step 3}. Formula (\ref{sub2}) holds. Denote by $S_m(a)=\{x\in\BBR^N:\abs{x_j}<m\}$. If $\overline S_R(a)\subset \CS$, the function 
$$(x,t)\mapsto c_\ga(t-\gt)^{-\frac{1+\ga}{q-1}}+w_{B_R}(x-a,t-\gt)
$$
is a supersolution of (\ref{1}) in $S_R(a)\ti (\gt,\infty)$ which is infinite on 
$ S_R(a)\ti\{\gt\}\cup\prt S_R(a)\ti [\gt,\infty)$ by \rprop{est2},
while $u$ is finite, thus it dominates $u$ in $S_R(a)\ti (\gt,\infty)$. Letting $\gt\to 0$ yields to
\bel {A1}
u(x,t)\leq c_\ga t^{-\frac{1+\ga}{q-1}}+w_{B_R}(x-a,t)\qquad\forall (x,t)\in Q_T^{S_R(a)}.
\ee
Conversely, the function 
$$(x,t)\mapsto u(x,t-\gt)+w_{B_R}(x-a,t-\gt)
$$
is a supersolution in $S_R(a)\ti (\gt,\infty)$ which dominates $c_\ga t^{-\frac{1+\ga}{q-1}}$ on $ S_R(a)\ti\{\gt\}\cup\prt S_R(a)\ti [\gt,\infty)$, thus as above, we obtain (\ref{A1}). Since 
$$\lim_{t\to 0}w_{B_R}(x-a,t)=0
$$
uniformly on $B_{R'}(a)$ for any $R'<R$, we derive (\ref{sub2}).
\qeda\medskip
\bprop{singset} For any relatively closed $\CS\in\Gw$, the set $\CU_{\CS,0}(\Gw)$ is not empty and it admits a minimal element  $\underline  u_{\CS,0}$ and a maximal element $\overline u_{\CS,0}$. 
\es
\Proof {\it Step 1: Existence of a maximal solution.} The maximal solution is constructed by thickening $\Gw$ and $\CS$ in defining for $0<\gs$
$$\BA {ll}
\Gw_\gs=\{x\in\BBR^N:\dist (x,\Gw)<\gs\},\qquad \CS_\gs=\{x\in\BBR^N:\dist (x,\overline \CS)\leq \gs\}
\EA$$
If $z\in \prt\Gw$, we denote by ${\bf n}_z$ the outward unit normal vector to $\Gw$ at $z$. Since $\prt\Gw$ is compact and $C^2$, there exists $\gs_0>0$ such that for any $(z,\gs)\in \prt\Gw\ti [0\gs_0]$, the mapping $\Gp:(z,\gs)\mapsto z+\gs{\bf n}_z$ is a $C^2$ diffeomorphism from $\prt\Gw\ti [0,\gs_0]$ to 
$\Gth'_{\gs_0}:=\overline \Gw_{\gs_0}\setminus\Gw$. The mapping $\Gp$ defines the flow coordinates near $\prt\Gw$. 

If $0<\gd<\gs$, there exists a unique solution $u=u_{n,\gs,\gd}$ of 
\bel{YZ0}\BA{ll}
\prt_tu-\Gd u+t^\ga u^q=0\qquad&\text {in }Q^{\Gw_\gs}_\infty\\
\phantom{\prt_t-\Gd u+t^\ga u^q}
u=0\qquad&\text {in }\prt \Gw_\gs\ti (0,\infty)\\
\phantom{\Gd u+t^\ga u^q}
u(.,0)=n\chi_{_{\CS_\gd}}  &\text {in } \Gw_\gs.
\EA\ee
Notice that $\CS_\gd$ is closed in $\Gw_\gs$ and $\inf\{\abs{z-z'}:z\in\CS_\gd,z'\in\Gw_\gd^c\}=\gd-\gs$. Existence is standard as well as uniqueness in the case where $\Gw$ is bounded. If $\Gw^c$ is bounded the proof goes as in the uniqueness proof in \rprop{existmea}. When $n\to\infty$, $\{u_{n,\gs,\gd}\}\uparrow u_{\gs,\gd}$ which is a solution of (\ref{1}) in $Q^{\Gw_\gs}_\infty$. Since 
$u_{n,\gs,\gd}$ satisfies (\ref{es0}), for any $r,\gt>0$ and any $a\in\prt\Gw_\gs$, $u_{n,\gs,\gd}$ remains uniformly continuous with respect to $n$ in $\overline Q^{\Gw_\gs}_\infty\cap \left(\overline B_r(a)\ti [2^{-2}\gt,2\gt]\right)$. Consequently $u_{\gs,\gd}(x,t)=0$ on $\prt\Gw_\gs\ti(0,\infty)\cap \left(\overline B_r(a)\ti [2^{-2}\gt,2\gt]\right)$. Therefore $u_{\gs,\gd}$ vanishes on $\prt\Gw_\gs\ti(0,\infty)$. If $a\in \Gw_\gs$ and $\dist (a,\CS_\gd)=r>0$, $u_{n,\gs,\gd}(x,t)\leq w_{B_r}(x-a,t)$ in $Q_{B_r(a)}^\infty$. This implies that $u_{n,\gs,\gd}$ remains uniformly continuous with respect to $n$ in $\overline B_{r'}(a)\ti [0,T)$ for any $0<r'<r$ and $T>0$. Since $u_{n,\gs,\gd}(x,t)\to 0$ in $\overline B_{r'}(a)$, $u_{\gs,\gd}$ inherits the same property. Consequently $u_{\gs,\gd}$ has initial trace $(\CS_\gd,0)$ in $\Gw_\gs$. By the maximum principle the mapping $(n,\gd)\mapsto u_{n,\gs,\gd}$ is increasing with respect to $n$ and decreasing with respect to $\gd$. Furthermore, if $0<\gd'<\gs'<\gs$ and $0<\gd<\gs$, there holds
$u_{n,\gs',\gd'}<u_{n,\gs,\gd}$ in $Q^{\Gw_{\gs'}}_\infty$, which implies $u_{\gs',\gd'}<u_{\gs,\gd}$. We define
\bel{YZ1}\BA{ll}
\displaystyle \overline  u_{\CS,0}=\lim_{\gs\to 0}\lim_{\gd\to 0}\lim_{n\to \infty}u_{n,\gs,\gd}.
\EA\ee
Then $\overline  u_{\CS,0}$ is a solution of (\ref{1}) in $Q^{\Gw}_\infty$. Since $\lim_{t\to 0}u_{\gs,\gd}(x,t)=0$ uniformly on any compact subset $K\subset\CR=\Gw\setminus\CS$, $\overline  u_{\CS,0}$ has initial trace $0$ on $\CR$. If $a\in \CS$, we denote by $\tilde u_{\infty,a}$ the function defined in $Q^{\Gw_\gs}_\infty$ by
$$
\tilde u_{\infty,a}(x,t)=\left\{\BA {ll}u_{\infty,a}(x,t)\qquad&\text{if }(x,t)\in Q^{\Gw}_\infty\\
0\qquad&\text{if }(x,t)\in Q^{\Gw_\gs}_\infty\setminus Q^{\Gw}_\infty
\EA\right.$$
Then $\tilde u_{\infty,a}$ is a subsolution of (\ref{1}) in $Q^{\Gw_\gs}_\infty$ which is smaller than $u_{\infty,a,\gs}$ which is the limit, when $k\to\infty$ of the solution $u_{k\gd_a,\gs}$ of
\bel{YZ2}\BA{ll}
\prt_tu-\Gd u+t^\ga u^q=0\qquad&\text {in }Q^{\Gw_\gs}_\infty\\
\phantom{\prt_t-\Gd u+t^\ga u^q}
u=0\qquad&\text {in }\prt \Gw_\gs\ti (0,\infty)\\
\phantom{\Gd u+t^\ga u^q}
u(.,0)=k\gd_a  &\text {in } \Gw_\gs.
\EA\ee
There holds, by \rth{sub},
$$u_{\gs,\gd}(x,t)\geq u_{\infty,a,\gs}\geq \tilde u_{\infty,a}(x,t)\qquad\forall (x,t)\in Q^{\Gw_\gs}_\infty.
$$
Letting successively $\gd\to 0$ and $\gs\to 0$ yields to $\overline  u_{\CS,0}\geq \tilde u_{\infty,a}=u_{\infty,a}$ in $Q^{\Gw_\gs}_\infty$. Therefore any $a\in\CS$ is a singular initial point of  $\overline  u_{\CS,0}$. Since $\CS\cup\CR=\Gw$, it follows that $tr_\Gw(u)=(\CS,0)$. Since $u_{\gs,\gd}$ satisfies (\ref{es0}) and $\prt\Gw_\gs$ has bounded curvature, independent of $\gs$, there holds classicaly
\bel{YZ3}\abs{\nabla u_{\gs,\gd}(x,t)}\leq ct^{-\frac{q+\ga}{q-1}}\qquad\forall (x,t)\in\prt \Gw_\gs\ti (0,\infty).
\ee
If $z\in \prt\Gw$ and then by the mean value theorem there exists $\gth\in (0,1)$ such that,
$$0\leq u_{\gs,\gd}(z,t)=u_{\gs,\gd}(z+\gs{\bf n}_z,t)-\gs\nabla u_{\gs,\gd}(z+\gth\gs{\bf n}_z,t).{\bf n}_z\leq c\gs t^{-\frac{q+\ga}{q-1}}.
$$
This implies that $\overline  u_{\CS,0}$ vanishes on $\prt\Gw_\gs\ti (0,\infty)$.

Let $u$ be any positive solution of (\ref{1}) in $Q^{\Gw}_\infty$, vanishing on $\prt \Gw\ti (0,\infty)$, with initial trace $(\CS,0)$. For $0<\gd<\gs$ fixed and for $R,\ge>0$, there exists $\gt_\ge>0$ such that, for any $\gt\in (0,\gt_\ge]$, 
$$u(x,\gt)\leq \ge\qquad\forall x\in \overline B_R\cap \Gw\setminus\CS_\gs. $$
This is due to the fact that $u(x,\gt)\to 0$ when $\gt\to 0$, uniformly on compact subset of $\overline B_R\cap \CR$. Assume that $\Gw$ is unbounded (the case where $\Gw$ is bounded is simpler since it does not require to introduce the barrier $w_{B_R}$) and let $R>0$ large enough so $\Gw^c\subset B_R$. By (\ref{es0})-(\ref{es1}) there exists $0<\gt_1\leq \gt_0$ such that for any $\gt\in (0,\gt_1]$, 
$$u(x,\gt)\leq w_{B_R}(x,\gt)+u_{\gs,\gd}(x,\gt)\qquad\forall x\in \overline B_R\cap \Gw\cap\CS_\gs.
$$
Furthermore $u(x,t)<w_{B_R}(x,t)$ for all $t>0$ and $x\in \prt B_R\cap \Gw$. Since $\ge+w_{B_R}+u_{\gs,\gd}$ is a supersolution for (\ref{1}) in $B_R\cap \Gw\ti (0,\infty)$, it follows that 
$$u(x,t)\leq \ge+w_{B_R}(x,t)+u_{\gs,\gd}(x,t)\qquad\forall (x,t)\in B_R\cap \Gw\ti (0,\infty).
$$
Letting successively $\gd\to 0$, $\gs\to 0$, $R\to \infty$ (here we use the fact that $w_{B_R}(x,t)\to 0$ when $R\to\infty$ by \rprop{est3}) and $\ge\to 0$ yields to $u\leq \overline  u_{\CS,0}$.
 \smallskip

\noindent {\it Step 2: Existence of a minimal solution.} The set $\CU_{\CS,0}(\Gw)$ is not empty since it contains $\overline  u_{\CS,0}$ and we may define
\bel{YZ4}
\tilde  u_{\CS,0}=\sup\{u_{\infty,a}:a\in\CS\},
\ee
and 
\bel{YZ5}
\hat  u_{\CS,0}=\inf\{u:u\in \CU_{\CS,0}(\Gw)\}.
\ee
The functions $\tilde  u_{\CS,0}$ and $\hat  u_{\CS,0}$ are respectively positive sub and supersolutions of (\ref{1}) in $Q^\Gw_\infty$. They are bounded from above by $\overline  u_{\CS,0}$ and from below by $u_{\infty,a}$ for any $a\in\CS$. Since $u_{\infty,a}\leq u$ for any  $a\in\CS$ and $u\in \CU_{\CS,0}(\Gw)$, it follows that $\tilde  u_{\CS,0}\leq \hat  u_{\CS,0}$. Therefore there exists a solution $\underline  u_{\CS,0}$ of (\ref{1}) in $Q^\Gw_\infty$ which satisfies 
and 
\bel{YZ6}
\tilde  u_{\CS,0}\leq \underline  u_{\CS,0}\leq \hat  u_{\CS,0}.
\ee
This implies that $\underline  u_{\CS,0}$ has initial trace $(\CS,0)$, it vanishes on $\prt\Gw\ti (0,\infty)$ and it is therefore the minimal element of $\CU_{\CS,0}(\Gw)$.\qeda\medskip

\noindent\Remark If $\dist (\CS,\Gw^c)>0$, it is not needed to replace $\Gw$ by a larger set $\Gw_\gs$ in order to construct the maximal solution. The construction of $\overline  u_{\CS,0}$ can be done in replacing $u_{n,\gs,\gd}$ by the solution  $u=u_{n,\gs}$ of 
\bel{YZ00}\BA{ll}
\prt_tu-\Gd u+t^\ga u^q=0\qquad&\text {in }Q^{\Gw}_\infty\\
\phantom{\prt_t-\Gd u+t^\ga u^q}
u=0\qquad&\text {in }\prt \Gw\ti (0,\infty)\\
\phantom{\Gd u+t^\ga u^q}
u(.,0)=n\chi_{_{\CS_\gd}}  &\text {in } \Gw,
\EA\ee
with $\gd<\gd_0:=\dist (\CS,\Gw^c)$.
\medskip

The next result is an extension of \rprop{stab}.
\bprop{stab2} Assume $\ga>-1$ and $1<q<q_{c,\ga}$. Let $\{u_n\}$ be a sequence of positive solutions of (\ref{1}) which converges to $u$ locally uniformly in $Q_T^\Gw$, and denote by $(\CS_n,\gm_n)$ and $(\CS,\gm)$ the respective initial trace of $u_n$ and $u$. If $\CA$ is an open subset of $\cap_n\CR_n$  and $\gm_n(\CA)$ remains bounded independently of $n\in \BBN$ (where $\CR_n=\Gw\setminus\CS_n$ and $\CR=\Gw\setminus\CS$), then $\CA\subset\CR$ and $\chi_{_\CA}\gm_n\rightharpoonup \chi_{_\CA}\gm$ in the weak sense of measures. Conversely, if $\CA\subset\CR$ , then  for any compact $K\subset \CA$, there exist $C_K>0$ and $n_K\in\BBN$ such that $\gm_n(K)\leq C_K$ for any $n\geq n_K$.
\es
\Proof Clearly (\ref{Ap0}) holds. We keep the notations of the proof of \rprop{stab}  where the first statement has been proved in assuming $\overline B_{\tilde r}(z)\subset\CA$. Since $\gm_n(\CA)$ remains bounded, there exists a subsequence $\{n_j\}$ and a positive measure $\gm'$ on $\CA$ such that $\gm_{n_j}\rightharpoonup \gm'$ in the weak sense of measures in $\CA$. Then $u_{\chi_{_{B_{\tilde r}(z)}} \gm_{n_j}}$ converges locally uniformly in $Q^{B_{\tilde r}(z)}_\infty$  to the solution 
$u_{\chi_{_{B_{\tilde r}(z)}} \gm'}$ of
\bel{Z0}\BA{ll}
\prt_tu-\Gd u+t^\ga u^q=0\qquad&\text {in }Q^{B_{\tilde r}(z)}_\infty\\
\phantom{\prt_t-\Gd u+t^\ga u^q}
u=0\qquad&\text {in }\prt B_{\tilde r}(z)\ti (0,\infty)\\
\phantom{\Gd u+t^\ga u^q}
u(.,0)=\chi_{_{B_{\tilde r}(z)}} \gm' &\text {in } B_{\tilde r}(z).
\EA\ee
Since $q<q_{c,\ga}$, the convergence of $u_{\chi_{_{B_{\tilde r}(z)}} \gm_{n_j}}$  and  $t^\ga u_{\chi_{_{B_{\tilde r}(z)}} \gm_{n_j}}^q$ respectively to $u_{\chi_{_{B_{\tilde r}(z)}} \gm_{n_j}}$ and $t^\ga u^q_{\chi_{_{B_{\tilde r}(z)}} \gm'}$ holds in $L^1(Q^{B_{\tilde r}(z)}_T)$ for any $T>0$. Relation (\ref{Ap7}) reads
\bel{Z2}
u_{\chi_{_{B_{\tilde r}(z)}} \gm_{n}}(x,t)\leq u_n(x,t)\leq u_{\chi_{_{B_{\tilde r}(z)}} \gm_{n}}(x,t)+w_{B_{\tilde r}(z)}(x,t)\qquad\text{in }Q^{B_{\tilde r}(z)}_\infty.
\ee
(see \rprop {est2}). Then $u_{n_j}$ and $t^\ga u^q_{n_j}$ converge to $u$ and $t^\ga u^q$ respectively,  in $L^1(Q^{B_{ r}(z)}_T)$ for any $r<\tilde r$. From (\ref{Ap3}), we derive 
\bel{Z3}
\dint_{\!\!Q^{B_{ r}(z)}_\infty}\left(-u(\gz_t+\Gd\gz)+\gz t^\ga u^q\right) dx dt=\myint{B_{ r}(z)}{}\gz(x,0) d\gm'(x),
\ee
for any $\gz\in C_c^{1,1;1}(\overline Q^{B_{ r}(z)}_\infty)$ which vanishes for $t$ large enough. This implies that $\gm'$ is the initial trace of $u$ in $B_{ r}(z)$, i.e. $\chi_{_{B_{ r}(z)}}\gm'=\chi_{_{B_{ r}(z)}}\gm$ and $\chi_{_{B_{ r}(z)}}\gm_n\rightharpoonup \chi_{_{B_{ r}(z)}}\gm$. Using a partition of unity, we conclude that $\chi_{_\CA}\gm_n\rightharpoonup \chi_{_\CA}\gm$.\smallskip

Conversely, we assume that there exist a compact set $K\subset \CA$ and a subsequence 
$\gm_{n_j}$ such that $\gm_{n_j}(K)\to\infty$. Thus, using the diagonal process, there exist $z\in K$ and another subsequence that we still denote $\gm_{n_j}$ such that 
$$\lim_{{n_j}\to\infty}\gm_{n_j}(B_\ge(z))=\infty\qquad\forall\ge>0.
$$
Therefore, we can construct a subsequence $\{n_{j_\ell}\}\subset \{n_{j}\}$ such that 
$$\gm_{n_{j_\ell}}(B_{2^{-n_{j_\ell}}}(z))=m_{n_{j_\ell}}\to\infty
$$
when $n_{j_\ell}\to\infty$. Since  the solution $u_{\chi_{_{B_{2^{-n_{j_\ell}}}(z)}} \gm_{n_{j_\ell}}}$ of
\bel{Z4}\BA{ll}
\prt_tu-\Gd u+t^\ga u^q=0\qquad&\text {in }Q^{B_{\tilde r}(z)}_\infty\\
\phantom{\prt_t-\Gd u+t^\ga u^q}
u=0\qquad&\text {in }\prt B_{\tilde r}(z)\ti (0,\infty)\\
\phantom{\Gd u+t^\ga u^q}
u(.,0)=\chi_{_{B_{2^{-n_{j_\ell}}}(z)}} \gm_{n_{j_\ell}} &\text {in } B_{\tilde r}(z).
\EA\ee
converges to $u^{B_{\tilde r}(z)}_{\infty\gd_z}$ which is the limit of the solution $u_{k\gd_z}$ of 
(\ref{Z4}) with initial data $u(.,0)=k\gd_z$, and is dominated by $u_{n_{j_\ell}}$ in $Q_T^{B_{\tilde r}(z)}$ we conclude that $u\geq u^{B_{\tilde r}(z)}_{\infty\gd_z}$ in $Q_T^{B_{\tilde r}(z)}$, which implies that $z\in\CS$, contradiction.\qeda
\bprop{local}
Assume $u_1$ and $u_2$ are two positive solutions of (\ref{1}) in $Q_\infty^\Gw$ with initial trace $(\CS,\gm)$. Then for any $a\in\CR$ and $R>0$ such that $B_R(a)\subset\CR$, there holds
\begin{equation}\label{loc1}\BA {ll}
\abs {u_1(x,t)-u_2(x,t)}\leq w_{B_R}(x-a,t)\qquad\forall (x,t)\in Q_\infty^{B_R(a)}
\EA\end{equation}
In particular $\lim_{t\to 0}\abs {u_1(x,t)-u_2(x,t)}=0$ uniformly on any compact subset of $\CR$.
\es
\Proof Since $u$ and $u'$ are solution of (\ref{1}) and $B_R(a)\in \CR$, for any $i=1,2$, $R'<R$ and $T>0$, 
$$\dint_{Q_T^{B_R'(a)}}t^\ga u^q_i(x,t)dx dt+\dint_{Q_T^{B_R'(a)}}u_i(x,t)dx dt<\infty,
$$
furthermore 
$$\lim_{t\to 0}\myint{B_R(a)}{}u_i(x,t)\gz(x) dx=\myint{B_R(a)}{}\gz(x)d\gm(x)\qquad\forall\gz\in C_c(B_R(a)).
$$
This implies that $u_i$ has a Sobolev trace on $\prt_\ell Q_T^{B_R'(a)}$ which belongs to $L^1$ and they are the limit, when $k\to\infty$ of the solutions $u_{i,k}$ of
\begin{equation}\label{loc2}\BA {ll}
\prt_tu-\Gd u+t^\ga \abs{u}^{q-1}u=0\qquad&\text{in }Q_\infty^{B_R'(a)}\\
\phantom{\prt_tu-\Gd u+t^\ga \abs{u}^{q-1}}
u=\min\{k,u_i\lfloor_{Q_\infty^{B_R'(a)}}\}\qquad&\text{in }\prt_\ell Q_\infty^{B_R'(a)}\\
\phantom{,\Gd u+t^\ga \abs{u}^{q-1}}
u(.,0)=\gm\qquad&\text{in }B_R'(a).
\EA\end{equation}
Since $u_{2,k}+w_{B_R'}(.-a)$ is a supersolution
$$u_{1,k}\leq u_{2,k}+w_{B_R'}(.-a)\Longrightarrow \abs{u_{1,k}-u_{2,k}}\leq w_{B_R'}(.-a)\text{ in }Q_\infty^{B_R'(a)}.$$
 Letting $k\to\infty$, $R'$ to $R$, we derive (\ref{loc1}). The second statement is a consequence of the fact that $\lim_{t\to 0}w_{B_R}(.-a)=0$, uniformly on $B_R'$ by \rprop{est2}.\qeda\medskip
 
 \noindent\Remark The previous estimate does not use the fact that $\prt\Gw$ is smooth and bounded. If the $u_i$ belong to $\CU_{\CS,\gm}(\Gw)$, estimate (\ref{loc1}) can be improved since the $u_i$ vanish on $\prt_\ell Q_\infty^{B_R(a)}$, and we obtain,
  \begin{equation}\label{loc3}\BA {ll}
\abs {u_1(x,t)-u_2(x,t)}\leq \min\left\{w_{B_R}(x-a,t),c_\ga t^{-\frac{1+\ga}{q-1}}\right\}\qquad\forall (x,t)\in Q_\infty^{B_R(a)}.
\EA\end{equation}

\bprop{existmea}Assume $\Gw\subseteq\BBR^N$ is either $\BBR^N$ or an open domain with a $C^2$ compact boundary, $\ga>-1$ and $1<q<q_{c,\ga}$. Then for any measure $\gm$ in $\Gw$ such that 
$\gm\lfloor_{\Gw_R}\in \mathfrak M_+^{b,\gr}(\Gw_R)$ where $\Gw_R=\Gw\cap B_R$, there exists a unique solution $u_\gm$ to
\begin{equation}\label{equmu}\BA {ll}
\prt_tu-\Gd u+t^\ga \abs{u}^{q-1}u=0\qquad&\text{in }Q^\Gw_\infty\\
\phantom{\prt_tu-\Gd u+t^\ga \abs{u}^{q-1}}
u=0\qquad&\text{in }\prt_\ell Q^\Gw_\infty\\
\phantom{,\Gd u+t^\ga \abs{u}^{q-1}}
u(.,0)=\gm\qquad&\text{in }\Gw,
\EA\end{equation}
and the mapping $\gm\mapsto u_\gm$ is increasing. Furthermore, if $\{\gm_n\}$ is a sequence of positive measures such that $\gm_n\lfloor_{\Gw_R}\in \mathfrak M_+^{b,\gr}(\Gw_R)$  which converges weakly to $\gm\lfloor_{\Gw_R}\in \mathfrak M_+^{b,\gr}(\Gw_R)$, then $\{u_{\gm_n}\}\to u_\gm$ locally uniformly in $Q^\Gw_\infty$.
\es
\Proof We recall that $u$ is a solution of (\ref{equmu}) if $u\in L^1_{loc}(\overline{Q^\Gw_\infty})$,  $\abs{u}^q\in L^1_{loc}(\overline{Q^\Gw_\infty};t^\ga \gr dxdt)$ satisfies
\begin{equation}\label{equmu1}\BA {ll}
\dint_{\!\!\!Q^\Gw_\infty}\left(-u(\prt_t\gz+\Gd\gz)+t^\ga \abs{u}^{q-1}u\gz\right)dx dt
=\myint{\Gw}{}\gz(.,0)d\gm
\EA\end{equation}
for any test function $\gz\in C^{2,1;1}_0(\overline{Q^\Gw_\infty})$. When $\Gw$ is bounded,  existence, uniqueness and stability are proved in  \cite{MV2}. Thus we assume that $\Gw$ is unbounded and we assume $R\geq R_0$ such that $\Gw^c\subset B_{R_0}$.  There exists a unique  solution $u_R$ of 
\begin{equation}\label{equmu2}\BA {ll}
\prt_tu-\Gd u+t^\ga \abs{u}^{q-1}u=0\qquad&\text{in }Q^{\Gw_R}_\infty\\
\phantom{\prt_tu-\Gd u+t^\ga \abs{u}^{q-1}}
u=0\qquad&\text{in }\prt_\ell Q^{\Gw_R}_\infty\\
\phantom{,\Gd u+t^\ga \abs{u}^{q-1}}
u(.,0)=\gm\lfloor_{\Gw_R}\qquad&\text{in }\Gw_R.
\EA\end{equation}
The function $u_R$ is nonnegative,  $R\mapsto u_R$ is increasing. For $R>R_1$, $u_R$ admits a Sobolev trace $f_{R_1}$ on $\prt B_{R_1}\ti (0,T)$ which is an integrable function, and $u_R$ is the unique solution of 
\begin{equation}\label{equmu3}\BA {ll}
\prt_tu-\Gd u+t^\ga \abs{u}^{q-1}u=0\qquad&\text{in }Q^{\Gw_{R_1}}_\infty\\
\phantom{\prt_tu-\Gd u+t^\ga \abs{u}^{q-1}}
u=0\qquad&\text{in }\prt \Gw\ti (0,\infty)\\
\phantom{\prt_tu-\Gd u+t^\ga \abs{u}^{q-1}}
u=f_{R_1}\qquad&\text{in }\prt B_{R_1}\ti (0,\infty)
\\
\phantom{,\Gd u+t^\ga \abs{u}^{q-1}}
u(.,0)=\gm\lfloor_{\Gw_{R_1}}\qquad&\text{in }\Gw_{R_1}.
\EA\end{equation}
Furthermore, $u_R\lfloor_{Q^{\Gw_{R_1}}_\infty}=\lim_{m\to\infty}v_m$, where $v_m$ is the unique solution of (\ref{equmu}) where the boundary data on $\prt B_{R_1}\ti (0,\infty)$ is replaced by $f_{R_1,m}=f_{R_1}\wedge m$ ($m\in\BBN^*$). Let  $v_{R_1}$ be the unique solution of 
\begin{equation}\label{equmu4}\BA {ll}
\prt_tu-\Gd u+t^\ga \abs{u}^{q-1}u=0\qquad&\text{in }Q^{\Gw_{R_1}}_\infty\\
\phantom{\prt_tu-\Gd u+t^\ga \abs{u}^{q-1}}
u=0\qquad&\text{in }\prt \Gw\ti (0,\infty)\\
\phantom{\prt_tu-\Gd u+t^\ga \abs{u}^{q-1}}
u=0\qquad&\text{in }\prt B_{R_1}\ti (0,\infty)
\\
\phantom{,\Gd u+t^\ga \abs{u}^{q-1}}
u(.,0)=\gm\lfloor_{\Gw_{R_1}}\qquad&\text{in }\Gw_{R_1},
\EA\end{equation}
If $w_{B_{R_1}}$ is the barrier function in $Q^{\Gw_{R_1}}_\infty$ which has been constructed in \rprop{est2}, $v_{R_1}+w_{B_{R_1}}$ is a supersolution for problem (\ref{equmu3}). Since it is larger than $v_m$ in $Q^{\Gw_{R_1}}_\infty$ for any $m>0$, there holds $u_R\leq v_{R_1}+w_{B_{R_1}}$, for any $R>R_1$. Then $u_R\uparrow u_\gm$ which is a solution of \ref{1} in $Q^{\Gw_{R}}_\infty$. By \rprop{est3}, $w_{B_{R_1}}$ remains uniformly bounded in $Q^{\Gw_{R'}}_\infty$ for any $R_0<R'<R_1$. Therefore $u_\gm$ shares the same property. If $\gz\in C_c^{1,1;1}(Q^{\overline \Gw}_\infty)$ vanishes on $Q^{\prt\Gw}_\infty)$ and for $\abs x>R'>R_0$, there holds for $R>R'>R_0$ and $T>0$,
\begin{equation}\label{equmu5}\BA {ll}
\dint_{Q^{\Gw}_\infty}\left(-u_R(x,t)(\prt_t\gz+\Gd\gz)+\gz t^\ga u_R^q\right)dx dt=\myint{\Gw}{}\gz(x,0)d\gm(x)-\myint{\Gw}{}\gz(x,T)u_R(x,T)dx
\EA\end{equation}
If we let $R\to\infty$ we deduce by the monotone convergence theorem that $u_\gm$ is a weak solution of (\ref{equmu}). This proves existence. 

For uniqueness, we consider $u_\gm$ and $u'_\gm$  two solutions of (\ref{equmu}). By the same argument as in the existence part, for any $R>0$, $u_\gm$ is smaller than the supersolution $u'_\gm+w_{B_{R_1}}$ in $Q^{\Gw_{R}}_\infty$. Since $\lim_{R\to\infty}w_{B_{R}}=0$ by \rprop {est3}  we obtain $u_\gm\leq u'_\gm$.
Similarly $u'_\gm\leq u_\gm$. Uniqueness implies the monotonicity of the mapping $\gm\mapsto u_\gm$.

For proving the stability, assume $\{\gm_n\lfloor_{\Gw_R}\}$ converges to $\gm\lfloor_{\Gw_R}$ in the weak sense of measures in $\mathfrak M_+^{b,\gr}(\Gw_R)$ for any $R>R_0$. Then the sequence of solutions $v_{n,R}$ of 
\begin{equation}\label{equmu6}\BA {ll}
\prt_tu-\Gd u+t^\ga \abs{u}^{q-1}u=0\qquad&\text{in }Q^{\Gw_{R}}_\infty\\
\phantom{\prt_tu-\Gd u+t^\ga \abs{u}^{q-1}}
u=0\qquad&\text{in }\prt \Gw\ti (0,\infty)\\
\phantom{\prt_tu-\Gd u+t^\ga \abs{u}^{q-1}}
u=0\qquad&\text{in }\prt B_{R}\ti (0,\infty)
\\
\phantom{,\Gd u+t^\ga \abs{u}^{q-1}}
u(.,0)=\gm_n\lfloor_{\Gw_{R}}\qquad&\text{in }\Gw_{R},
\EA\end{equation}
converges to the solution $v_R$ of (\ref{equmu4}) with $R_1=R$. In particular $v_{n,R}\to v_R$ and $t^\ga v^q_{n,R}\to t^\ga v^q_{R}$ in $L^1(Q^{\Gw_{R}}_T)$ and by standard regularity result the convergence of $v_{n,R}$ towards $v_R$ holds uniformly on $\overline \Gw_{R}\ti[\ge,T]$ for any $0<\ge<T$. Furthermore $u_{\gm_n}\leq v_{n,R}+W_{B_{R}}$ in $Q^{\Gw_{R}}_\infty$. This jointly with standard local regularity results for heat equation, implies that $\{u_{\gm_n}\}$ remains uniformly bounded and hence relatively compact for the topology of uniform convergence on any compact set of $\overline \Gw\ti[\ge,\infty)$. Thus there exist a subsequence $\{u_{\gm_{n_k}}\}$ and a function $u^*\in C^(\overline Q^\Gw_\infty)$ such that $u_{\gm_{n_k}}\to u^*$ locally uniformly in $\overline\Gw\ti(0,\infty)$. Since $t^\ga u^q_{\gm_n}\leq t^\ga v^q_{n,R}+t^\ga W^q_{B_{R}}$, there also holds by the dominated convergence theorem $t^\ga u^q_{\gm_{n_k}}\to t^\ga {u^*}^q$ in $L^1_{loc}(\overline\Gw\ti[0,\infty))$. Henceforth letting  $n_k\to\infty$ in the expression
\begin{equation}\label{equmu7}\BA {ll}
\dint_{Q^{\Gw}_\infty}\left(-u_{\gm_{n_k}}(x,t)(\prt_t\gz+\Gd\gz)+\gz t^\ga u^q_{\gm_{n_k}}\right)dx dt\\[4mm]\phantom{-----------}=\myint{\Gw}{}\gz(x,0)d\gm_{n_k}(x)-\myint{\Gw}{}\gz(x,T)u_{n_k}(x,T)dx,
\EA\end{equation}
 where $\gz\in C_c^{1,1;1}(Q^{\overline \Gw}_\infty)$, we conclude that $u^*=u_\gm$ and that $u_{\gm_{n}}\to u_\gm$.\qeda\medskip

\bprop{meas} Assume $F$ is a non-empty relatively closed subset of $\Gw$, $\CR=\Gw\setminus F$ and $\gm\in \mathfrak M_+(\CR)$. If we set
\begin{equation}\label{M1}\BA {ll}
\prt_\gm F=\{z\in F:\gm(\CR\cap B_r(z))=\infty,\,\forall r>0\}, 
\EA\end{equation}
then $\prt_\gm F$ is relatively closed in $\Gw$. If $\CR^\gm=\Gw\setminus \prt_\gm F$, it contains $\CR$ and if $\gm^*$ is the measure defined in  $\CR^\gm$ by $\gm$ on $\CR$ and $0$ in $\CR^\gm\cap \CR^c$, then 
there exist a minimal positive solution $\underline u_{\gm^*}$ and a maximal solution $\overline u_{\gm^*}$ of (\ref{1}) vanishing on $\prt\Gw\ti (0,\infty)$  satisfying $tr_\Gw(u)=(\prt_\gm F,\gm^*)$.
Furthermore $\underline u_{\gm^*}$ and $\overline u_{\gm^*}$ are respectively the minimal and the maximal element of $\CU_{\prt_\gm F,\gm^*}(\Gw)$.
\es
\Proof The set $\prt_\gm F$ is the blow-up set of the measure $\gm$. It is clearly a relatively closed subset of $\Gw$ included into $\overline\CR\setminus\CR$. \smallskip

\noindent{\it Step 1: Existence of a minimal solution.} For $\gd>0$, we denote $(\prt_\gm F)_\gd:=\{x\in\Gw:\dist (x,\prt_\gm F)\leq \gd\}$ and $\CR_\gd^\gm=\Gw\setminus(\prt_\gm F)_\gd\subset \CR^\gm$. We define the Radon measure $\gm_\gd$ on  $\Gw$ by 
$$
\gm_\gd=\left\{\BA {ll}\gm\qquad&\text {on}\; \CR_\gd^\gm\cap\CR\\
0\qquad&\text {on}\;  F\cup (\prt_\gm F)_\gd\EA\right.
$$ 
Then $\gm_\gd$ is a positive Radon measure in $\Gw$ and by \rprop{existmea} problem (\ref{equmu}) with initial data $\gm_\gd$ admits a unique positive solution $u_{\gm_\gd}$. Furthermore the mapping $\gd\mapsto u_{\gm_\gd}$ is nonincreasing, and we set $\underline u_{\gm^*}=\lim_{\gd\to 0}u_{\gm_\gd}$. Then $\underline u_{\gm^*}$ is a positive solution of (\ref{1}) in $Q^{\Gw}_\infty$ which vanishes on $\prt\Gw\ti (0,\infty)$ and has initial trace $(\CS',\gm')$. If $a\in\CR^\gm$, there exists $R>0$ such that $\overline B_R(a)\subset \CR^\gm$ and $\gd_a>0$ such that $\overline B_R(a)\subset \CR^\gm_\gd$ for $0<\gd<\gd_a$, that we assume below. By the maximum principle there holds
$$v_{\gm_\gd}(x,t)\leq u_\gd(x,t)\leq v_{\gm_\gd}(x,t)+w_{B_R}(x-a,t)\qquad\text{in }Q_\infty^{B_R(a)},
$$
where $v_{\gm_\gd}$ is the solution of 
$$\BA {ll}
\prt_tu-\Gd u+t^\ga u^q=0\qquad&\text{in }Q_\infty^{B_R(a)}\\
\phantom{\prt_t-\Gd u+t^\ga u^q}
u\geq 0\qquad&\text{in }Q_\infty^{B_R(a)}\\
\phantom{\Gd u+t^\ga u^q}
u(x,t)=0\qquad&\text{in }\prt B_R(a)\ti (0,\infty)\\
\phantom{ -u+t^\ga u^q}
u(.,0)=\gm_\gd\chi{_{B_R(a)}}\qquad&\text{in }B_R(a).
\EA$$
Letting $\gd\to 0$, then $\gm_\gd\chi{_{B_R(a)}}\uparrow \gm^*\chi{_{B_R(a)}}$, which yields to
$$v_{\gm^*\lfloor_{B_R(a)}}(x,t)\leq \underline u_{\gm^*}(x,t)\leq v_{\gm^*\lfloor_{B_R(a)}}(x,t)+w_{B_R}(x-a,t)\qquad\text{in }Q_\infty^{B_R(a)}.
$$
By a partition of unity, it implies that for any $\gz\in C_c(\CR^\gm)$, we have
$$\lim_{t\to 0}\myint{\CR^\gm}{}\underline u_{\gm^*}(x,t)\gz(x) dx=\myint{\CR^\gm}{}\gz d\gm^*(x).
$$
Therefore $\CR^\gm\subset\CR'$ and $\gm'\lfloor_{\CR^\gm}=\gm^*$. If $z\in \CR'\cap (\CR^\gm)^c=\CR'\cap\prt_\gm F$, $\gm^* (B_r(z)\cap\CR^\gm) =\gm (B_r(z)\cap\CR)=\infty$ for any $r>0$ while there exists $r_0>0$ such that $\gm'(B_{r_0}(z))<\infty$. By \rprop{stab2}, for any $r'<r_0$ there exists $C>0$ such that $\gm_\gd(\overline B_r(z))\leq C$. By the monotone convergence theorem, it implies  $\gm(\overline B_r(z)\cap \CR)\leq C$, which contradicts the definition of $\prt_\gm F$. Thus 
$\CS'=\prt_\gm F$ and $tr_\Gw(\underline u_{\gm^*})=(\prt_\gm F,\gm^*)$.

Let us assume that $\Gw$ is unbounded, $R_0$ is such that $\Gw^c\subset B_{R_0}$ and $\Gw_R=\Gw\cap B_R$ for $R>R_0$. Let $\gd,\ge>0$, there exists $\gt_\ge$ such that 
\begin{equation}\label{M3}u_{\gm_\gd}(x,\gt)\leq\ge\qquad\forall (x,\gt)\in \Gw_R\cap (\prt_\gm F)_{\frac{\gd}{2}}\ti[0,\gt_\ge]
\end{equation}
Let $u\in\CU_{\prt_\gm F,\gm^*}(\Gw)$. In order to compare $u_{\gm_\gd}$ and $v:=u+w_{B_R}+\ge$ in $\Gw_R\setminus (\prt_\gm F)_{\frac{\gd}{2}}\ti(0,\gt_\ge]$ we see that 
$$u_{\gm_\gd}(.,0)=\chi_{_{\CR^\gm_R\setminus (\prt_\gm F)_{\gd} }}\gm^*\leq u(.,0)=\chi_{_{\CR^\gm_R\setminus (\prt_\gm F)_{\frac{\gd}{2}} }}\gm^*$$
and both are bounded Radon measures. Since $u_{\gm_\gd}(t,0)\leq v(x,t)$ in $\prt(\Gw_R\setminus (\prt_\gm F)_{\frac{\gd}{2}})\ti(0,\gt_\ge]$ and $v$ is a supersolution, it follows that  $u_{\gm_\gd}\leq v$ in $\Gw_R\setminus (\prt_\gm F)_{\frac{\gd}{2}}\ti(0,\gt_\ge]$. Using (\ref{M3}) we conclude that 
 \begin{equation}\label{M4'}
 u_{\gm_\gd}(x,\gt)\leq v(x,\gt)\qquad\forall (x,\gt)\in \Gw_R\ti[0,\gt_\ge].
\end{equation}
Next, applying the comparison principle in $\Gw_R\ti [\gt^*,\infty)$ between the solution $u_{\gm_\gd}$ and the supersolution $u+\ge+w_{B_R}$, we conclude that (\ref{M4'}) holds in $\Gw_R\ti [\gt,\infty)$
 and thus in $Q_\infty^{B_R}$. Letting successively $R\to\infty$, $\ge\to 0$ and $\gd\to 0$, we conclude that $\underline u_{\gm^*}\leq u$, thus $u_{\gm^*}$ is the minimal element of $\CU_{\prt_\gm F,\gm^*}(\Gw)$. 
\smallskip


\noindent{\it Step 2: Existence of a maximal solution.} Let $\gd>0$ and $u\in \CU_{\prt_\gm F,\gm^*}(\Gw)$. By \rprop{local}, for any $R>0$ and $\ge>0$ there exists $\gt_\ge$ such that 
\begin{equation}\label{M4}u(x,t)\leq \underline u_{\gm^*}(x,t)+\ge\qquad\forall (x,t)\in \Gw_R\setminus (\prt_\gm F)_\gd \ti (0,\gt_\ge],
\end{equation}
and by (\ref{es0}), $u(x,t)\leq c_\ga t^{-\frac{1+\ga}{q-1}}$. Let $\gt\in (0,\gt_\ge]$ and $w_{\gd,\gt}$ be the solution of (\ref{equmu}) in $Q^{\Gw_R}_\infty$ with initial data $\gm$ replaced by
\begin{equation}\label{M5}
h_{\gd,\gt}=\left\{\BA {ll} \underline u_{\gm^*}(x,\gt)\qquad&\text{if }x\in \Gw_R\setminus (\prt_\gm F)_\gd
\\ c_\ga\gt^{-\frac{1+\ga}{q-1}}&\text{if }x\in \Gw_R\cap (\prt_\gm F)_\gd
\EA\right.
\end{equation}
By (\ref{M4}), (\ref{M5}) and the maximum principle, 
\begin{equation}\label{M6}
u(x,t+\gt)\leq w_{\gd,\gt}(x,t)+\ge+w_{B_R}(x,t+\gt)\qquad\forall (x,t)\in Q^{\Gw_R}_\infty.
\end{equation}
Let $\overline u_{(\prt_\gm F)_\gd}$ be the maximal element of $\CU_{(\prt_\gm F)_\gd,0}(\Gw)$, which exists by \rprop{singset}. Then, by (\ref{es1}), for any $\gd'>\gd$, there exists $\gt_{\gd'}\in (0,\gt_\ge]$ such that for any $\gt\in (0,\gt_{\gd'}]$
\begin{equation}\label{M7}
\max\{\underline u_{\gm^*}(.,.+\gt),\overline u_{(\prt_\gm F)_\gd}(.,.+\gt)\}\leq w_{\gd,\gt}\leq \underline u_{\gm^*}(.,.+\gt)+\overline u_{(\prt_\gm F)_{\gd'}}(.,.+\gt)\qquad\text {in } Q^{\Gw_R}_\infty.
\end{equation}
Up to a sequence $\{\gt_n\}$ converging to $0$, $\{w_{\gd,\gt_n}\}$ converges, locally uniformly in $Q^\Gw_\infty$ to a solution $w_{\gd}$ of (\ref{1}) in $Q^\Gw_\infty$ which satisfies
\begin{equation}\label{M8}
\max\{\underline u_{\gm^*},\overline u_{(\prt_\gm F)_\gd}\}\leq w_{\gd}\leq \underline u_{\gm^*}+\overline u_{(\prt_\gm F)_{\gd'}}\qquad\text { in } Q^{\Gw_R}_\infty.
\end{equation}
We can replace $\gd'$ by $\gd$ in this inequality, this proves that  $w_{\gd}$ vanishes on $\prt_\ell Q^{\Gw_R}_\infty$ has initial trace $(\gm_\gd,(\prt_\gm F)_\gd)$, therefore (\ref{M6}) becomes
\begin{equation}\label{M9}
u(x,t)\leq w_{\gd}(x,t)+\ge+w_{B_R}(x,t)\qquad\forall (x,t)\in Q^{\Gw_R}_\infty.
\end{equation}
Letting successively $R\to\infty$ and $\ge\to 0$ we deduce that $w_{\gd}$ is larger that any $u\in \CU_{\prt_\gm F,\gm^*}(\Gw)$ in $Q^{\Gw_R}_\infty$. Since $h_{\gd,\gt}$ decreases with $\gd$, $w_{\gd}$ shares this property and the limit, denoted by $\overline u_{\gm^*}$ is a solution of (\ref{1}) in $Q^{\Gw}_\infty$ which vanishes on $\prt_\ell Q^{\Gw_R}_\infty$ which is large than $u$, thus it is the maximal element of $\CU_{\prt_\gm F,\gm^*}(\Gw)$.\qeda\medskip


\bprop {meas-alt} Under the assumptions of $\rprop{meas}$, we set $ F_\gd:=\{x\in\Gw:\dist(x, F)\leq\gd\}$ and $\CR_\gd:=\Gw\setminus F_\gd\subset\CR$. If we define the measure $\tilde\gm_\gd$ in $\Gw$ by
$$
\tilde\gm_\gd=\left\{\BA {ll}\gm\qquad&\text {on}\; \CR_\gd\\
0\qquad&\text {on}\;  F_\gd,\EA\right.
$$ 
then $u_{\tilde \gm_\gd}\uparrow\underline u^*$ when $\gd\downarrow 0$.
\es
\Proof There holds $\tilde\gm_\gd\leq \gm_\gd$ which implies $u_{\tilde\gm_\gd}\leq u_{\gm_\gd}$. When $\gd\to 0$, 
$u_{\tilde\gm_\gd}\uparrow\underline u^*\leq \underline u_{\gm^*}$, thus $\underline u^*$ is a positive solution of (\ref{1}) in $Q^\Gw_\infty$ which vanishes on $\prt_\ell Q^\Gw_\infty$. Then $tr_\Gw(\underline u^*)=(\CS'',\gm'')$ and 
$\CS''\subset \prt_\gm F$ and  $\gm''\leq \gm^*$ on $\CR^\gm$. Furthermore $\gm''=\gm=\gm^*$ on $\CR$, as in the proof of \rprop{meas}. Since $\gm^*=0$ on $\CR^\gm\setminus\CR$ it follows that $\gm''= \gm^*$ on $\CR^\gm$. Let $a\in \CR''\cap\prt_\gm F$ and $R>0$ such that $\overline B_R(a)\subset \CR''$. Then $\gm''(\overline B_R(a))<\infty$. Therefore
$$\gm''(\overline B_R(a)\cap\CR)=\gm(\overline B_R(a)\cap\CR)<\infty,
$$
contradiction. Thus $\CS''=\prt_\gm F$ and $tr_\Gw(\underline u^*)=(\prt_\gm F, \gm^*)$. Since $\underline u^*\leq \underline u_{\gm^*}$ and $\underline u_{\gm^*}$ is minimal, it follows that $\underline u^*= \underline u_{\gm^*}$.\qeda

\bprop{subsup}Assume $v\in C(\overline\Gw\ti (0,\infty))$ is a positive sub-solution (resp. supersolution) of (\ref{1}) in $Q^\Gw_\infty$ which vanishes on $\prt\ell Q^\Gw_\infty$. Then there exists a minimal solution $\gp_+( v)$ larger than $v$ (resp. a maximal solution $\gp_-( v)$ smaller than $v$ and vanishing on $\prt\ell Q^\Gw_\infty$).
\es
\Proof Assume $v$ is a subsolution. Let $\gt>0$ and let $u_\gt$ be the solution of
\bel{sub1+}\BA{ll}
\prt_tu-\Gd u+t^\ga u^q=0\qquad&\text{in }\Gw\ti(\gt,\infty)\\
\phantom{--,,--}
u\geq 0, u=0\qquad&\text{in }\prt\Gw\ti(\gt,\infty)\\
\phantom{-,----}
u(.,\gt)=v(.,\gt)\qquad&\text{in }\Gw.
\EA\ee
Existence and uniqueness follows from \rprop{existmea}. Furthermore $u_\gt\geq v$ in $\Gw\ti(\gt,\infty)$. This implies that for $0<\gt<\gt'$, $u_{\gt}\geq u_{\gt'}$. Since $u_{\gt}(x,t)\leq c_\ga(t-\gt)^{-\frac{1+\ga}{q-1}}$,  there exists $\gp_+( v)=\lim_{\gt\to 0}u_\gt$, and  $\gp_+( v)$ is a positive solution of  (\ref{1}) in $Q^\Gw_\infty$ and is larger than $v$. If $u$ is any positive solution of  (\ref{1}) in $Q^\Gw_\infty$,  vanishing on $\prt\ell Q^\Gw_\infty$ and larger than $v$, for any $\gt>0$ it is larger than $u(.,\gt)$, thus it is larger than $u_\gt$ on $\Gw\ti(\gt,\infty)$. Therefore $u\geq \gp_+( v)$.

Assume now that $u$ is a supersolution. We define $u_\gt$ by (\ref{sub1+}). Then $u_\gt\leq v$ and 
$u_{\gt}\leq u_{\gt'}$ for $0<\gt<\gt'$. Then $\gp_-( v)=\lim_{\gt\to 0}u_\gt$, and  $\gp_-( v)$ is a positive solution of  (\ref{1}) in $Q^\Gw_\infty$ and is smaller than $v$, and thus vanishing on $\prt\ell Q^\Gw_\infty$. Similarly as above $\gp_-( v)$ is larger than any positive solution smaller than $v$.\qeda 

\bth{exist} Assume $\Gw\subseteq\BBR^N$ is either $\BBR^N$ or an open domain with a $C^2$ compact boundary, $\ga>-1$ and $1<q<q_{c,\ga}$. Then for any $\gn\approx (\CS,\gm)\in  \frak B^{reg}(\Gw)$ there exist a maximal positive solution $\overline u_{\CS,\gm}$ and a minimal positive solution $\underline u_{\CS,\gm}$ of (\ref{1}) in $Q^\Gw_\infty$ vanishing on $\prt\Gw\ti (0,\infty)$ with initial trace $\gn$. Furthermore, if  $\inf\{\abs{z-z'}:z\in\CS,z'\in\Gw^c\}>0$, then
\bel{h8}
\overline u_{\CS,\gm}-\underline u_{\CS,\gm}\leq 
\overline u_{\CS,0}-\underline u_{\CS,0}.
\ee
\es 
\Proof {\it Step 1: Construction of the maximal and minimal solutions.} The functions $\underline u_{\gm^*}$,  $\underline u_{\CS,0}$, $\overline u_{\gm^*}$  and $\overline u_{\CS,0}$ have been defined in \rprop{meas} and \rprop{singset}. Since $\sup\{\underline u_{\gm^*},\underline u_{\CS,0}\}$ is a subsolution of (\ref{1}) which is smaller than the supersolution $\overline u_{\gm^*}+\overline u_{\CS,0}$
we set
 \bel{h9}\BA{rrr}
(i)\phantom{-}\underline u_{\CS,\gm}=\gp_+(\sup\{\underline u_{\gm^*},\underline u_{\CS,0}\})\phantom{-}\text{and }\phantom{-}
(ii)\phantom{-}\overline u_{\CS,\gm}=\gp_-(\overline u_{\gm^*}+\overline u_{\CS,0}).
\EA\ee
Then $\underline u_{\CS,\gm}$ and $\overline u_{\CS,\gm}$ are solutions which satisfy
 \bel{h10}
 \sup\{\underline u_{\gm^*},\underline u_{\CS,0}\}\leq \underline u_{\CS,\gm}\leq \overline u_{\CS,\gm}
 \leq \overline u_{\gm^*}+\overline u_{\CS,0}.
 \ee
 Therefore $\underline u_{\CS,\gm}$ and $\overline u_{\CS,\gm}$ vanish on $\prt_\ell Q^\Gw_\infty$, they have initial trace $\gm$ on $\CR$ and are larger than any $u_{\infty,a}$ for $a\in\CS$ (notice that $\prt_\gm\CS\subset\CS$). This implies that $\underline u_{\CS,\gm}$ and $\overline u_{\CS,\gm}$ belong to $\CU_{\CS,\gm}(\Gw)$. 
 
 Let $u\in \CU_{\CS,\gm}(\Gw)$. Then for $\gs>\gd>0$ and $\ge,R>0$, there exists $\gt_1>0$ such that 
 $$u(x,t)\leq u_{\gs,\gd}(x,t) \text{ for }  (x,t)\in (\CS_{\frac{\gd}{2}}\cap\Gw_R)\ti (0,\gt_1].$$
  There exists $\gt_2\in (0,\gt_1]$ such that 
  $$u(x,t)\leq \overline u_{\gm^*}(x,t)+\ge+w_{B_R}(x,t) \text{ for }  (x,t)\in  (\Gw_R\setminus \CS_{\frac{\gd}{2}})\ti (0,\gt_2].$$
  Therefore
  $$u(x,t)\leq u_{\gs,\gd}(x,t)+\overline u_{\gm^*}(x,t)+\ge+w_{B_R}(x,t) \text{ for }  (x,t)\in  \Gw_R\ti (0,\gt_2].$$
  Therefore 
   \bel{h11}
   u\leq u_{\gs,\gd}+u_{\gm^*}+\ge+w_{B_R}\quad\text{ in }  Q^{\Gw_R}_\infty
    \ee
    Letting successively $R\to\infty$, $\ge\to 0$, $\gd\to 0$ and $\gs\to 0$ we obtain $  u\leq \overline u_{\CS,0}+\overline u_{\gm^*}$, and therefore $u\leq \gp_+(\overline u_{\CS,0}+\overline u_{\gm^*})=\overline u_{\CS,\gm}$. Next, we also have
    $$u\geq \tilde u_{\CS,0}:=\sup\{u_{\infty,a}:a\in \CS\}
    \Longrightarrow u\geq \underline u_{\CS,0}=\gp_+(\tilde u_{\CS,0}),
    $$
  by (\ref{YZ4}). With the notations of \rprop{meas-alt} with $F=\CS$, for any $R>0$, $\gd>0$ and $\ge>0$, there exists  $\gt_\ge$ such that 
  $$u_{\tilde \gm_\gd}(x,t)\leq u(x,t)+w_{B_R}(x,t)+\ge\qquad\text{in }\Gw_R\ti (0,\gt_\ge],
  $$
  because the support of $\tilde\gm_\gd$ is included in $\Gw\setminus \CS_\gd$. Therefore this last inequality holds in $Q^{\Gw_R}_\infty$ and consequently 
  $$\sup\{u_{\tilde\gm_\gd},\underline u_{\CS,0}\}\leq u+\ge+w_{B_R}\qquad\text{in }Q^{\Gw_R}_\infty,
  $$
  and we can let $R\to\infty$ and $\ge\to 0$ to obtain $\sup\{u_{\tilde\gm_\gd},\underline u_{\CS,0}\}\leq u$ in $Q^{\Gw}_\infty$.
  Letting $\gd\to 0$ and using \rprop{meas-alt} we get $\sup\{\underline u_{\gm^*},\underline u_{\CS,0}\}\leq \underline u_{\CS,\gm}\leq u$.\medskip
  
  
    \noindent  {\it Step 2: Alternative construction.} For $0<\gd<\gs$ and $n\in\BBN_*$ we denote by $\overline u_{n,\gs,\gd,\gm}$ the solution of
   \bel{W1}\BA {ll}
\prt_tu-\Gd u+t^\ga u^q=0\qquad&\text{in }Q_\infty^{\Gw_\gs}\\
\phantom{\prt_t-\Gd u+t^\ga u^q}
u\geq 0\qquad&\text{in }Q_\infty^{\Gw_\gs}\\
\phantom{\Gd u+t^\ga u^q}
u(x,t)=0\qquad&\text{in }\prt_\ell Q_\infty^{\Gw_\gs}\\
\phantom{ -u+t^\ga u^q}
u(.,0)=\gm_\gd+n\chi{_{\CS_\gd}}\qquad&\text{in }\Gw_\gs.
\EA\ee
We denote here  $\Gw_\gs=\{x\in \BBR^N:\dist(x,\Gw)<\gs\}$, $\CS_\gd=\{x\in \BBR^N:\dist(x,\CS)<\gd\}$ and 
$\CR_\gd=\Gw\cap \CS^c_\gd$ and $\gm_\gd=\chi{_{\CR_\gd}}\gm$. The same arguments of monotonicity as in \rprop{singset} and \rprop{meas} show that 
     \bel{W2}\BA {ll}\displaystyle
\lim_{\gs\to 0}\lim_{\gd\to 0}\lim_{n\to \infty}=\overline u_{n,\gs,\gd,\gm}=\overline u_{\CS,\gm}.
\EA\ee
If $\gt>0$ we denote by $\underline u_{\gt,\gd,\gm}$ the solution of
   \bel{W3}\BA {ll}
\prt_tu-\Gd u+t^\ga u^q=0\qquad&\text{in }Q_\infty^{\Gw}\\
\phantom{\prt_t-\Gd u+t^\ga u^q}
u\geq 0\qquad&\text{in }Q_\infty^{\Gw}\\
\phantom{\Gd u+t^\ga u^q}
u(x,t)=0\qquad&\text{in }\prt_\ell Q_\infty^{\Gw}\\
\phantom{ -u+t^\ga u^q}
u(.,0)=\gm_{\gd}+\chi_{_{\CS_\gd}}\underline u_{\CS,0}(.,\gt) \qquad&\text{in }\Gw,
\EA\ee
Using estimate (\ref{sub1}) and \rprop{local} is is easy to prove that $\underline u_{\gt,\gd,\gm}\leq u$ for any $u\in\CU_{\CS,0}(\Gw)$. Furthermore 
$$\max\{u_{\gm_\gd},u_{\chi_{_{\CS_\gd}}\underline u_{\CS,0}(.,\gt)}\}\leq \underline u_{\gt,\gd,\gm}\leq u_{\gm_\gd}+\underline u_{\CS,0}(.,.+\gt),$$
since have
 $$u_{\chi_{_{\CS_\gd}}\underline u_{\CS,0}(.,\gt)}\leq \underline u_{\CS,0}(.,\gt)\leq 
u_{\chi_{_{\CS_\gd}}\underline u_{\CS,0}(.,\gt)}+C\left(\frac{\gd^2}{\gt}\right)^{\frac{1+\ga}{q-1}-\frac{1}{2}}e^{-\frac{\gd^2}{4\gt}}
 $$
 by (\ref{sing4}) with $N=1$. Set $c(\gd,\gt)=C\left(\frac{\gd^2}{\gt}\right)^{\frac{1+\ga}{q-1}-\frac{1}{2}}e^{-\frac{\gd^2}{4\gt}}$, then  
$$ \max\{u_{\gm_\gd},u_{\chi_{_{\CS_\gd}}\underline u_{\CS,0}(.,\gt)}\}
\geq \max\{u_{\gm_\gd},\underline u_{\CS,0}(.,\gt)-c(\gd,\gt)\}\geq \max\{u_{\gm_\gd},\underline u_{\CS,0}(.,\gt)\}-c(\gd,\gt).
$$
Therefore, if $\underline u_{\gt_n,\gd,\gm}\to \underline u_{0,\gd,\gm}$ locally uniformly in $Q_\infty^{\Gw}$, then $\underline u_{0,\gd,\gm}$ is a solution of (\ref{1}) in $Q_\infty^{\Gw}$ which satisfies
   \bel{W3+1}
   \max\{u_{\gm_\gd},\underline u_{\CS,0}\}\leq \underline u_{\gd,0,\gm}\leq u_{\gm_\gd}+\underline u_{\CS,0},
\ee
and is smaller than any $u\in\CU_{\CS,0}(\Gw)$. There exists $\gd_n\to 0$ such that $\underline u_{0,\gd_n,\gm}\to \underline u_{0,0,\gm}$. Then 
   \bel{W3+2}
   \max\{\underline u_{\gm^*},\underline u_{\CS,0}\}\leq \underline u_{0,0,\gm}\leq \underline u_{\gm^*}+\underline u_{\CS,0},
\ee
and $\underline u_{0,0,\gm}$ is an element of $\CU_{\CS,0}(\Gw)$ smaller than any $u\in\CU_{\CS,0}(\Gw)$. Thus $\underline u_{0,0,\gm}=\underline u_{\CS,\gm}$ and 
   \bel{W3+3}
  \lim_{\gd\to 0}  \lim_{\gt\to 0}\underline u_{\gt,\gd,\gm}=\underline u_{\CS,\gm}.
  \ee
\medskip


  \noindent  {\it Step 3: Proof of (\ref{h8}).} We assume $\inf\{\abs{z-z'}:z\in\CS,z'\in\Gw^c\}=\gd_0>0$, so that we can take $\gs=0$ in the construction of $\overline u_{\CS,\gm}$. Put $\gt=(\frac{c_\ga}{n})^{\frac{q-1}{1+\ga}}$ and 
\bel{W5}\overline Z_{n,\gd,\gm}=\overline u_{n,\gd,\gm}-\overline u_{n,\gd,0}\,,\;
\underline Z_{\gt,\gd,\gm}=\underline u_{\gt,\gd,\gm}-\underline u_{\gt,\gd,0}\;\text{and }\;W_{n,\gd,\gm}=\overline Z_{n,\gd,\gm}- \underline Z_{\gt,\gd,\gm}.
  \ee
Then $w=W_{n,\gd,\gm}$ satisfies
\bel{W6}\BA {l}\prt_tw-\Gd w+t^{\ga}\left(\overline u^q_{n,\gd,\gm}-\overline u^q_{n,\gd,0}-\underline u^q_{\gt,\gd,\gm}+\underline u^q_{\gt,\gd,0}
\right)
\EA\ee
in $Q^\Gw_\infty$ and we can write
\bel{W7}\BA {l}
\overline u^q_{n,\gd,\gm}-\overline u^q_{n,\gd,0}-\underline u^q_{\gt,\gd,\gm}+\underline u^q_{\gt,\gd,0}=(\overline u^q_{n,\gd,\gm}-\underline u^q_{\gt,\gd,\gm})-(\overline u^q_{n,\gd,0}-\underline u^q_{\gt,\gd,0})\\\phantom{\overline u^q_{n,\gd,\gm}-\overline u^q_{n,\gd,0}-\underline u^q_{\gt,\gd,\gm}+\underline u^q_{\gt,\gd,0}}
=d_\gm(\overline u_{n,\gd,\gm}-\underline u_{\gt,\gd,\gm})-d_0(\overline u_{n,\gd,0}-\underline u_{\gt,\gd,0})
\EA\ee
where 
\bel{W8}d_\gm(x,t)=\left\{\BA {ll}
\myfrac{\overline u^q_{n,\gd,\gm}-\underline u^q_{\gt,\gd,\gm}}{\overline u_{n,\gd,\gm}-\underline u_{\gt,\gd,\gm}}\quad&\text{if }\overline u_{n,\gd,\gm}\neq \underline u_{\gt,\gd,\gm}\\[4mm]
0&\text{if }\overline u_{n,\gd,\gm}= \underline u_{\gt,\gd,\gm}
\EA\right.\ee
and $d_0$ is defined accordingly. Since 
$$\overline u_{n,\gd,\gm}\geq \max\{ \underline u_{\gt,\gd,\gm},\overline u_{n,\gd,0}\}\;\text{and }\;  \underline u_{\gt,\gd,0}\leq \min\{\underline u_{\gt,\gd,\gm},\overline u_{n,\gd,0}\},
$$
there holds $d_\gm\geq d_0\geq $ by convexity. Using the fact that $\overline u_{n,\gd,0}-\underline u_{\gt,\gd,0} \geq 0$ is infers that 
$$d_\gm(\overline u_{n,\gd,\gm}-\underline u_{\gt,\gd,\gm})-d_0(\overline u_{n,\gd,0}-\underline u_{\gt,\gd,0})\geq 
d_\gm(\overline u_{n,\gd,\gm}-\underline u_{\gt,\gd,\gm}-\overline u_{n,\gd,0}+\underline u_{\gt,\gd,0})
$$
Finally (\ref{W6}) becomes
\bel{W9}
\prt_tw-\Gd w+t^{\ga}d_\gm w\leq 0\qquad\text{in }Q^\Gw_\infty.
\ee
Furthermore, in the sense of measures, 
$$w(.,0)=\gm_\gd+n\chi{_{\CS_\gd}}-n\chi{_{\CS_\gd}}-(\gm_{\gd}+\chi_{_{\CS_\gd}}\underline u_{\CS,0}(.,\gt))+\chi_{_{\CS_\gd}}\underline u_{\CS,0}(.,\gt)=0.
$$
Because $w=0$ in $\prt_\ell Q^\Gw_\infty$ it follows $w\leq 0$ by the maximum principe. Therefore
\bel{W10}
\overline Z_{n,\gd,\gm}\leq 
\underline Z_{\gt,\gd,\gm}\Longrightarrow
\overline u_{n,\gd,\gm}-\underline u_{\gt,\gd,\gm}\leq\overline u_{n,\gd,0}-\underline u_{\gt,\gd,0}.
\ee
If we let successively $n\to\infty$ (and therefore $\gt\to 0$) and $\gd\to 0$, we obtain (\ref{h8}).\qeda\medskip


\noindent\Remark We do not know if (\ref{h8}) holds if we do not assume $\inf\{\abs{z-z'}:z\in\CS,z'\in\Gw^c\}=\gd_0>0$. However, if for $\gth>0$ we set $\CS_\gth=\CS\cap\{x\in\Gw:\dist (x,\Gw^c)\geq\gth\}$, then we have
\bel{W11}
\overline u_{\CS_\gth,\gm}-\underline u_{\CS_\gth,\gm}\leq 
\overline u_{\CS_\gth,0}-\underline u_{\CS_\gth,0}.
\ee
Furthermore all the four above functions increases when $\gth$ decreases to 0. If we set
\bel{W12}\BA {lll}
&(i)\qquad\qquad&\lim_{\gth\to 0}\overline u_{\CS_\gth,\gm}=\overline u_{\underline \CS,\gm}\\
&(ii)\qquad\qquad&\lim_{\gth\to 0}\overline u_{\CS_\gth,0}=\overline u_{\underline \CS,0}\\
&(iii)\qquad\qquad&\lim_{\gth\to 0}\underline u_{\CS_\gth,\gm}=\underline u_{\underline \CS,\gm}\\
&(iv)\qquad\qquad&\lim_{\gth\to 0}\underline u_{\CS_\gth,0}=\underline u_{\underline \CS,0},
\EA
\ee
then we infer that
\bel{W13}
\overline u_{\underline \CS,\gm}-\underline u_{\underline \CS,\gm}\leq 
\overline u_{\underline \CS,0}-\underline u_{\underline \CS,0}.
\ee\qeda

Our final result is the following existence and uniqueness theorem the proof is close to the one of \cite[Th 3.5]{MV1}, therefore we present only the main ideas and the needed changes.
\bth{exist-uniq} Assume $\Gw\subseteq\BBR^N$ is either $\BBR^N$ or an open domain with a $C^2$ compact boundary, $\ga>-1$ and $1<q<q_{c,\ga}$. Then for any $\gn\approx (\CS,\gm)\in  \frak B^{reg}(\Gw)$ such that $\gr\gm$ is bounded in any neighborhood of $\prt\Gw$ and $\inf\{\abs{z-z'}:z\in\CS,z'\in\Gw^c\}>0$, the set $\CU_{\CS,\gm}(\Gw)$ contains one and only one element.
\es
\Proof  We assume that $\CS\neq\{\emptyset\}$ otherwhile uniqueness is already known. Thanks to (\ref{h8}) it is enough to prove that $\overline u_{ \CS,0}=\underline u_{\CS,0}$. \smallskip

\noindent {\it Case 1: $\Gw=\BBR^N$}. For  $\gs>0$ and $a\in\CS$ set 
$\CP_\gs(a)=\{(x,t)\in Q_\infty:\abs{x-a}\leq \gs\sqrt t\}$ and $\CP_\gs=\cup_a\in\CP_\gs(a)$. Because of (\ref{es0}), (\ref{sing4}), (\ref{sub1}) and (\ref{init4}), for any $\gs>0$, there exists $C_\gs>1$ such that
\bel{Un1}
\overline u_{\CS,0}(x,t)\leq C_\gs\underline u_{\CS,0}(x,t)\qquad\forall (x,t)\in \CP_\gs.
\ee
Fix $\gs>0$. If $y\in \BBR^N\setminus\CS$, we set $r(y)=\dist (y,\CS)$. Using \rprop{est-1} and estimate (\ref{ss4}), we obtain that 
\bel{Un2}
\overline u_{\CS,0}(x,t)\leq 2Nt^{-\frac{1+\ga}{q-1}}W_\ga\left(\myfrac{r(y)-|x-y|}{\sqrt t}\right)\qquad
\forall (x,t)\in Q_\infty^{B_{r(y)}(y)}.
\ee
If we take $x=y$ in the above estimate, we derive from (\ref{ss3}), 
\bel{Un3}
\overline u_{\CS,0}(y,t)\leq C_1\left(\myfrac{r(y)}{\sqrt t}\right)^{2\frac{1+\ga}{q-1}}e^{-\frac{r^2(y)}{4 t}}(1+o(1)).
\ee
Since there exists $z\in \CS\cap\prt B_{r(y)}(y)$, we have also from (\ref{sub1}) and (\ref{sing4})
\bel{Un4}
\underline u_{\CS,0}(y,t)\geq C_2\left(\myfrac{r(y)}{\sqrt t}\right)^{2\frac{1+\ga}{q-1}-N}e^{-\frac{r^2(y)}{4 t}}(1+o(1)).
\ee
Thus there exists $C_3>1$ such that 
\bel{Un5}
\overline u_{\CS,0}(y,t)\leq \underline u_{\CS,0}(y,C_3t).
\ee
For $\gt>0$, we denote by $u_{1,\gt}$ and $u_{2,\gt}$ the solutions of (\ref{1}) in $Q_{\gt,\infty}:=\BBR^N\ti (\gt,\infty)$ with respective initial data
\bel{Un6}\BA {ll}
u_1(x,\gt)=C_\gs\underline u_{\CS,0}(x,\gt)\chi_{_{\CP_\gs}}(x,\gt)
\\[2mm]u_2(x,\gt)=\underline u_{\CS,0}(x,C_3\gt)(1-\chi_{_{\CP_\gs}}(x,\gt)).
\EA\ee
Then $u_1+u_2$ is a supersolution of (\ref{1}) in $Q_{\gt,\infty}$. Therefore
\bel{Un7}
(u_1+u_2)(x,\gt)\geq \overline u_{\CS,0}(x,\gt)\;\forall x\in \BBR^N\Longrightarrow
(u_1+u_2)(x,t)\geq \overline u_{\CS,0}(x,t)\;\forall (x,t)\in Q_{\gt,\infty}.
\ee
Since $C_\gs\underline u_{\CS,0}$ is a supersolution of (\ref{1}) which is larger than $u_1$ at $t=\gt$, then $u_1\leq C_\gs\underline u_{\CS,0}$ in $Q_{\gt,\infty}$. Next the function
$(x,t)\mapsto w(x,t):=\underline u_{\CS,0}(x,t+(C_3-1)\gt)$ satisfies
$$\BA {ll}\prt_tw-\Gd w+h(t,\gt)t^\ga w^q=0\quad&\text{in }Q_{\gt,\infty}\\
[2mm]\phantom{;;-\Gd w==,=\gt)}
w(.,\gt)=\underline u_{\CS,0}(.,C_3\gt)\quad&\text{in }\BBR^N
\EA$$
where 
$$h(t,\gt)=\left(\myfrac{t+(C_3-1)\gt}{t}\right)^\ga.
$$
Then
$$h(t,\gt)\leq C_4:=\max\{1,C_3^\ga\}\qquad\forall t\geq\gt.
$$
This implies that $C^{1/(q-1)}_4w$ is a supersolution of (\ref{1}) in $Q_{\gt,\infty}$ which is larger than 
$u_2$ for $t=\gt$. Thus it dominates $u_2$ in $Q_{\gt,\infty}$. Combining the above estimates on $u_1$ and $u_2$ with (\ref{Un7}), we derive that for any $\gt>0$, there holds
\bel{Un8}
\overline u_{\CS,0}(x,t)\leq C_5\left(\underline u_{\CS,0}(x,t+(C_3-1)\gt)+\underline u_{\CS,0}(x,t)\right)
\qquad\forall (x,t)\in Q_{\gt,\infty}.
\ee
where $C_5=\max\left\{C^{1/(q-1)}_4,C_\gs\right\}$. Letting $\gt\to 0$ we finally obtain the key estimate
with $C=2C_5$,
\bel{Un9}
\overline u_{\CS,0}\leq C\underline u_{\CS,0}
\quad\text{ in } Q_{\infty}.
\ee 
Then end of the proof is the same as in \cite[Th 3.5]{MV1}, but we recall it for the sake of completeness: If $\overline u_{\CS,0}\neq\underline u_{\CS,0}$, then strict inequality holds. For $0<\gb<C^{-1}$ the function 
$\tilde u=\underline u_{\CS,0}-\ga (\overline u_{\CS,0}-\underline u_{\CS,0})$ is a  supersolution of (\ref{1}) in  $Q_{\infty}$ (this due to the convexity of $r\mapsto r^q$) which satisfies 
$\gb\underline u_{\CS,0}\leq \tilde u<\underline u_{\CS,0}$. For $0<\gamma<\gb$, the function $\gamma\underline u_{\CS,0}$ is a subsolution smaller than $\tilde u$. Then there exists a solution $u'$ of  (\ref{1}) in  $Q_{\infty}$ which satisfies
\bel{Un10}
\gamma\underline u_{\CS,0}\leq u'\leq \tilde u<\underline u_{\CS,0}.
\ee 
This implies that $tr_{\BBR^N}(u')=(\CS,0)$, which contradicts the minimality of $\underline u_{\CS,0}$.\smallskip

\noindent{\it Case 2: $\prt\Gw$ is nonempty and compact.} Again the proof is similar to the one of \cite[Th 3.5]{MV1}. We denote by $\underline u^\Gw_{\CS,0}$ and $\overline u^\Gw_{\CS,0}$ the minimal and the maximal solutions of (\ref{1}) in  $Q^\Gw_{\infty}$ with initial trace $(\CS,0)$ and we assume that 
$\dist (\CS,\Gw^c)=\gd>0$. We also set $\underline u^{\BBR^N}_{\CS,0}=\underline u_{\CS,0}$ and $\overline u^{\BBR^N}_{\CS,0}=\overline u_{\CS,0}$. Clearly $\overline u^\Gw_{\CS,0}\leq \overline u_{\CS,0}$ in $Q^\Gw_{\infty}$. Furthermore, if we denote by $k(t)$ the maximum of $\overline u_{\CS,0}(x,t)$ for $x\in\prt\Gw$, then $\lim_{t\to 0}k(t)=0$ (this is due to the fact that $\dist (\CS,\Gw^c)>0$). Clearly, the construction of the minimal solutions shows that $\underline u_{\CS,0}(x,t)\leq \underline u^\Gw_{\CS,0}(x,t)+k(t)$ in $Q^\Gw_{\infty}$. Therefore
\bel{Un11}
\overline u^\Gw_{\CS,0}(x,t)\leq \underline u^\Gw_{\CS,0}(x,t)+k(t)\quad\forall (x,t)\in Q^\Gw_{\infty}.
\ee 
If we fix $\gt>0$, $t\mapsto \underline u^\Gw_{\CS,0}(x,t)+k(t)$ is a supersolution of (\ref{1}) in $Q^\Gw_{\infty}$ which is larger that $\overline u^\Gw_{\CS,0}(x,t)$ at $t=\gt$. This implies 
\bel{Un12}
\overline u^\Gw_{\CS,0}(x,t)\leq \underline u^\Gw_{\CS,0}(x,t)+k(\gt)\quad\forall (x,t)\in Q^\Gw_{\gt,\infty}.
\ee 
If we let $\gt\to 0$, we deduce that $\overline u^\Gw_{\CS,0}=\underline u^\Gw_{\CS,0}$.\qeda

\end{document}